\documentclass[12pt,thmsa,titlepage]{article}
\usepackage{amsmath}
\usepackage{amstext}
\usepackage{amsfonts}
\usepackage{amssymb}
\usepackage{graphicx}
\newtheorem{theorem}{Theorem}

\newtheorem{lemma}[theorem]{Lemma}

\newtheorem{proposition}[theorem]{Proposition}
\newtheorem{remark}[theorem]{Remark}

\newenvironment{proof}[1][Proof]{\textbf{#1.} }{\ \rule{0.5em}{0.5em}}
\begin{document}

\author{M. L. Frankel\thanks{ Department of Mathematical Sciences, Indiana
University--Purdue University Indianapolis, Indianapolis, IN 45205} and V.
Roytburd\thanks{Department of Mathematical Sciences, Rensselaer Polytechnic
Institute, Troy, NY 12180-3590 \textbf{Corresponding author, a LaTeX file of
the article is available from this author}; \textbf{tel. (518)-276-6889; email
roytbv@rpi.edu}}}
\title{Finite-Dimensional Attractor for a Nonequilibrium Stefan Problem with Heat Losses}
\date{\empty}
\maketitle
\begin{abstract}
We study a two-phase\ modified Stefan problem modeling solid combustion and
nonequilibrium phase transition. The problem is known to exhibit a variety of
non-trivial dynamical scenarios. We develop a priori estimates and establish
well-posedness of the problem in weighted spaces of continuous functions. The
estimates secure sufficient decay of solutions that allows for an analysis in
Hilbert spaces. We demonstrate existence of compact attractors in the weighted
spaces and prove that the attractor consists of sufficiently regular
functions. This allows us to show that the Hausdorff dimension of the
attractor is finite.



\textbf{Running head}: FINITE-DIMENSIONAL ATTRACTOR FOR A NONEQUILIBRIUM
STEFAN PROBLEM
\end{abstract}

%To be submitted to J Differential Equations, finished on March 23, 2003
%For JDE:elsart.cls (use this file if you are using LaTeX2e, the current version of LaTeX), 
%elsart.sty and elsart13.sty (use these two files if you are using LaTeX3.09, the previous version of LaTeX), 
%guidelines for users of elsart, a template file for quick start, and the instruction booklet 
%"Preparing articles with LaTeX." 
%

\section{Introduction}

The subject of this paper is a study of dynamics of a nonequilibrium two-phase
Stefan problem modeling condensed phase combustion and some phase transition
processes. It was demonstrated numerically \cite{port2} that the sharp
interface model of the condensed phase combustion also known as
Self-propagating High-temperature Synthesis (SHS) generates a remarkable
variety of complex thermokinetic oscillations. In addition to its theoretical
interest SHS finds technological applications as a method of synthesizing
certain technologically advanced materials, see \cite{munir}, \cite{var2} and
also \cite{var1} for a popular exposition. The process is characterized by
highly exothermic reactions propagating through mixtures of fine elemental
reactant powders, resulting in the synthesis of compounds. The dynamical
scenarios exhibited by the model include a Hopf bifurcation, period doubling
cascades leading to chaotic pulsations, a Shilnikov-Hopf bifurcation etc.
These scenarios are well-known for the finite-dimensional dynamical systems
and suggest a possibility that the essential dynamics of the free-interface
problem may be finite-dimensional as well. Indeed, we \ have been able to
prove \cite{comp,sell} that compactness and finite dimensionality of the
attractor\ take place for a simpler one-phase problem.

However, the methods of the papers dealing with the one-phase problem are not
directly applicable to the sharp interface problem of the condensed phase
combustion which is the subject of the present paper. The principal difficulty
that arises here as compared to the one-phase problem is that the presence of
the additional \ temperature field behind the propagating interface (in the
product phase) creates an additional degree of freedom that is not easily
controllable. This difficulty is overcome in the present paper; \textit{we
show that Hausdorff dimension of the attractor is finite}. The paper draws on
the approach of our previous work \cite{comp2} discussing compactness of the
attractor; in addition to the dimension estimate, results presented in this
paper clarify structure and regularity properties of the attractor.

There is a substantial literature that treats analytical aspects of the
initial--boundary value problem for different sharp-interface models with
kinetics \ related to the problem (\ref{he}-\ref{jc}) below, see \cite{luc,
radkevich, chen,yin1}. These works are concerned with basic issues of mostly
local in time existence. We also note recent papers by Brauner \textit{et al.
and }Lorenzi\textit{, }\cite{brauner,brauner1, lorenzi}, which study
weakly-nonlinear dynamical behavior of solutions of related problems. In
particular they consider perturbations of traveling-wave initial data and
investigate their instability and bifurcations. In contrast, the principal
focus of the present paper is in strongly nonlinear asymptotic dynamics for a
wide range of initial data and parametric regimes.

The free-interface problem is formulated as follows: find $s(t)$ and $u(x,t)$
such that
\begin{gather}
u_{t}=u_{xx}-\gamma u,\quad x\neq s(t),\quad t>0,\label{he}\\
u(x,0)=u_{0}(x)\geq0,\label{ic}\\
g[u(s(t),t)]=v(t),\label{kc}\\
\lbrack u_{x}(s(t),t)]:=u_{x}^{+}(s(t),t)-u_{x}^{-}(s(t),t)=v(t)\label{jc}%
\end{gather}
where $v(t)$ is the interface velocity, $s(t)=\int_{0}^{t}v(\tau)d\tau$ is its
position, $u$ is the temperature, and the derivatives $u_{x}^{+}$ and
$u_{x}^{-}$ are taken from right side and left side of the free interface
respectively. The last term in the heat equation (\ref{he}) is due to the heat
losses into the medium surrounding the combustible or solidifying substance
via Newton's cooling law with a non-dimensional coefficient $\gamma>0$.

Dynamics of the physical system is determined by the feedback mechanism
between the heat release due to the kinetics $g(u|_{x=s(t)})$\ and the heat
dissipation by the medium. The second interface condition (\ref{jc}) (the
Stefan boundary condition) expresses the balance between the heat produced at
the free boundary and its diffusion by the adjacent medium. As the problem
describes, generally speaking, propagation of a phase transition front, the
first interface condition (\ref{kc}) is a manifestation of the
\textit{nonequilibrium} nature of the transition; its analog for the classical
Stefan problem is just $u|_{x=s(t)}=0.$\ We should mention that in contrast
with the nonequilibrium problem, the dynamics of the classical Stefan problem
is relatively trivial. The surrounding matter is assumed to be at the
temperature of the fresh combustible mixture at $-\infty$ (the original phase
in the phase transition interpretation). By the same token the heat loss will
reduce the temperature in the product phase to that of the medium. Thus the
behavior of the solution at infinity should satisfy $\lim
\nolimits_{x\rightarrow\pm\infty}u(x,t)=0$.

In order to estimate the Hausdorff dimension of the attractor we need to
develop some additional technical tools. We first develop a priori estimates
and establish well-posedness of the problem in weighted spaces of continuous
functions. These estimates, which constitute the analytical core of the paper,
secure sufficient decay of solutions that allows us to carry out analysis in a
Hilbert space. It should be noted that volume estimates, which are the basis
for the Hausdorff dimension bound, require a Hilbert structure in the
underlying space. Next we are able to extend our results on existence of
compact attractors \cite{comp2} to the weighted spaces. After that we study
the evolution on the attractor and prove that the semigroup on the attractor
is onto and one-to-one: it yields that the attractor consists of sufficiently
regular functions. As a consequence we are able to demonstrate that problem is
well-posed and its attractor is precompact in a Hilbert space.

This allows us to estimate the Hausdorff dimension of the attractor based on
the techniques described for instance in \cite{temam}. We study evolution of
the infinitesimal volume along the trajectories in the attractor and
demonstrate that for sufficiently large $m$ that is defined solely by the
physical properties of the problem, the $m$-dimensional volume decays
exponentially. After that we prove that the semigroup is uniformly
differentiable which, combined with the estimate for the linearized evolution
of the infinitesimal volume leads to the conclusion that the Hausdorff
dimension of the attractor is finite.

\section{Properties of solutions: Previous results}

In this section we present some pertinent background information from
\cite{minqu} (certain statements are slightly modified and clarified). The
following theorem summarizes existence results:

\begin{theorem}
Suppose that the kinetic functions $g$ satisfies the following
assumptions:\newline $(A1)\ $ $g(u)$ is a continuously differentiable,
monotone decreasing, negative function on $(0,\infty)$ with $g(0)=-v_{0}$ for
some velocity $-v_{0}<0;$\newline $(A2)$ $\ g(u)$ is sublinear: ${\lim
}_{u\rightarrow\infty}g(u)/u=0;$ \newline and that the initial data
$u_{0}(x)\in C(-\infty,\infty)$. Then there exists one and only one classical
solution of the free interface problem (\ref{he})-(\ref{jc}). The solution is
uniformly bounded for all $t>0$.
\end{theorem}

The proof is based on the reduction to an integral equation for the interface
velocity
\begin{equation}
v(t)=g\left(  e^{-\gamma t}\int_{-\infty}^{\infty}G(s(t),t,\xi,0)u_{0}%
(\xi)d\xi-\int_{0}^{t}G(s(t),t,s(\tau),\tau)e^{-\gamma(t-\tau)}v(\tau
)d\tau\right)  ,\label{ie}%
\end{equation}
which arises from the interface condition and the single layer potential
representation for the solution operator:
\begin{equation}
(Tu_{0})(x,t):=u(x,t)=e^{-\gamma t}\int_{-\infty}^{\infty}G(x,t,\xi
,0)u_{0}(\xi)d\xi-\int_{0}^{t}G(x,t,s(\tau),\tau)e^{-\gamma(t-\tau)}%
v(\tau)d\tau,\label{iu}%
\end{equation}
where
\begin{equation}
G(x,t,\xi,\tau)=\exp\{-\frac{(x-\xi)^{2}}{4(t-\tau)}\}[4\pi(t-\tau)]^{-1/2}%
\end{equation}
is the heat kernel and $s(t)=\int_{0}^{t}v(\tau)d\tau$.

In the sequel we replace the sublinearity condition $(A2)$ by a stronger
condition. We assume that $g(u)$ is a monotonically decreasing differentiable
function on $[0,\infty]$ with $|g^{\prime}|\leq C$ and satisfying
\begin{equation}
-V^{0}\leq g(u)\leq-v_{0}\;\mathrm{{for}\;{some}\;}V^{0},v_{0}%
>0.\label{kinetics}%
\end{equation}
These conditions are satisfied, for instance, for the standard Arrhenius
kinetics where $v=V^{0}\exp(-A/(u-u_{\infty}))$.

Under some additional conditions the following smoothness result holds
\cite{minqu}:

\begin{theorem}
Let the initial data $u_{0}$ be twice differentiable in $x<0$ and $x>0$ with
bounded derivatives and satisfy the matching condition:
\begin{equation}
g(u_{0}(0))=\frac{\partial u_{0}^{+}}{\partial x}(0)-\frac{\partial u_{0}^{-}%
}{\partial x}(0),
\end{equation}
in addition, let the derivative of the kinetics function $g^{\prime}$ be
Lipschitz continuous. Then the velocity $v$ is differentiable.
\end{theorem}

\section{A priori estimates in weighted spaces}

For our purposes we need to establish certain a priori bounds on the solution
in appropriate weighted spaces that are introduced next. Let$\ \omega_{\alpha
}$ be\ the weight $\omega_{\alpha}(x)=e^{\alpha|x|}$; we define
\[
|\,f\,|_{\alpha}=\sup(\omega_{\alpha}(x)|f(x)|),\;C_{\alpha}=\{f\in
C(-\infty,\infty):\;|\,f\,|_{\alpha}<\infty\}
\]
We will demonstrate that the global existence results can be can be extended
to $C_{\alpha}$ (note the obvious imbedding, if $\beta>\alpha\geq0$ then
$C_{\beta}\subset C_{\alpha}\subset C(-\infty,\infty)$)$.$ Similarly we define
Hilbert versions of weighted spaces:
\[
\left\|  f\,\right\|  _{\alpha}=\left\|  \omega_{\alpha}f\,\right\|
_{L_{2}(-\infty,\infty)},\;H_{\alpha}=\{f:\;\left\|  f\,\right\|  _{\alpha
}<\infty\}
\]
It is easy to see that if $\beta<\alpha$ then $C_{\alpha}\subset H_{\beta},$
and for any $f\in C_{\alpha}$
\begin{equation}
\left\|  f\,\right\|  _{\beta}=\left(  \int_{-\infty}^{\infty}\omega_{\beta
}^{2}f^{2}dx\right)  ^{1/2}=\left(  \int_{-\infty}^{\infty}\frac{\omega
_{\beta}^{2}}{\omega_{\alpha}^{2}}\omega_{\alpha}^{2}f^{2}dx\right)
^{1/2}\leq|\,f\,|_{\alpha}\frac{1}{\sqrt{\alpha-\beta}}.\label{beta-alpha}%
\end{equation}

It is convenient to split the representation formula (\ref{iu}) for the
semigroup operator $T$ into two operators: the contribution of the free
boundary
\begin{equation}
T_{1}(t)u^{0}(x^{\prime})=-\int_{0}^{t}e^{-\gamma(t-\tau)}G(x^{\prime
},t,s(\tau),\tau)\left[  v(\tau)\right]  d\tau
\end{equation}
and that of the initial data
\begin{equation}
T_{2}(t)u^{0}(x^{\prime})=e^{-\gamma t}\int_{-\infty}^{\infty}G(x^{\prime
},t,\xi,0)u^{0}(\xi)d\xi
\end{equation}

In the sequel we will frequently encounter integrals of the error function
type. To estimate them we employ the following simple result.

\begin{lemma}
\label{errfun} For $a,b>0$
\[
\int_{a}^{\infty}\exp(-b\eta^{2})d\eta\leq\left\{
\begin{array}
[c]{l}%
\frac{1}{2\sqrt{b}}\exp(-ba^{2}),\text{ for }a>1/\sqrt{b}\\
\frac{\sqrt{\pi}}{2\sqrt{b}},\text{ for }0\leq a<1/\sqrt{b}%
\end{array}
\right.
\]
\end{lemma}

\begin{proof}
If $a\sqrt{b}>1$ then
\[
\int_{a}^{\infty}\exp(-b\eta^{2})d\eta=\frac{1}{\sqrt{b}}\int_{a\sqrt{b}%
}^{\infty}\exp(-\eta^{2})d\eta\leq\int_{a}^{\infty}\eta\exp(-b\eta^{2}%
\}d\eta=\frac{1}{2\sqrt{b}}\exp(-ba^{2})
\]
On the other hand
\[
\int_{a}^{\infty}\exp(-b\eta^{2})d\eta\leq\int_{0}^{\infty}\exp(-b\eta
^{2})d\eta=\frac{\sqrt{\pi}}{2\sqrt{b}}%
\]
\end{proof}

\subsection{Estimates for the solution: Contribution from initial data}

\begin{proposition}
For sufficiently small $\alpha$ (if $\alpha$ satisfies $\alpha^{2}+\alpha
V^{0}-\gamma<0$) the contribution from the initial data in the $C_{\alpha}%
$-norm decays exponentially in time:
\[
|u_{2}(.,t)|_{\alpha}\leq2\exp[(-\gamma+\alpha^{2}+\alpha V^{0})t]\,|u_{0}%
|_{\alpha}%
\]
\end{proposition}

\begin{proof}
For the contribution from the initial data, $u_{2}(.,t)=T_{2}(t)u_{0},$ we
have:
\begin{align*}
& |u_{2}(.,t)|_{a}=\sup_{x}\left\{  \omega_{\alpha}(x-s(t))\;|\int_{-\infty
}^{\infty}e^{-\gamma t}G(x,t,\xi,0)u_{0}(\xi)d\xi|\right\} \\
& =\sup\left\{  e^{-\gamma t}\omega_{\alpha}(x-s(t))|\int_{-\infty}^{\infty
}\frac{1}{\omega_{\alpha}(\xi)}G(x,t,\xi,0)u_{0}(\xi)\omega_{\alpha}(\xi
)d\xi|\right\} \\
& \leq\frac{e^{-\gamma t}}{2\sqrt{t\pi}}|u_{0}|_{\alpha}\sup\{\omega_{\alpha
}(x-s(t))[\int_{-\infty}^{0}e^{\alpha\xi}\exp(-\frac{(x-\xi)^{2}}{4t}%
)d\xi+\int_{0}^{\infty}e^{-\alpha\xi}\exp(-\frac{(x-\xi)^{2}}{4t})d\xi]\}
\end{align*}
Each of the integrals should be estimated separately for $x<0$ and $x>0$ with
the maximum of the estimates chosen for the estimate of the norm$.$ At the
same time it is easy to see that the two integrals can be transformed into
each other through the change $x\rightarrow-x,$ therefore it suffices to
estimate only one of them and double the result.

We proceed as follows
\[
\frac{1}{2\sqrt{t}}\int_{0}^{\infty}e^{-\alpha\xi}\exp(-\frac{(x-\xi)^{2}}%
{4t})d\xi=e^{\alpha^{2}t}e^{-\alpha x}\int_{-\frac{(x-2t\alpha)}{2\sqrt{t}}%
}^{\infty}\exp(-\eta^{2}\}d\eta
\]
For $x>0$ the above expression
\[
e^{\alpha^{2}t}e^{-\alpha x}\int_{-\frac{(x-2t\alpha)}{2\sqrt{t}}}^{\infty
}\exp(-\eta^{2}\}d\eta\leq\sqrt{\pi}e^{\alpha^{2}t}e^{-\alpha x}%
\]
\ and therefore
\begin{gather*}
\frac{e^{-\gamma t}}{\sqrt{\pi}}|u_{0}|_{\alpha}\sup_{x>0}\omega_{\alpha
}(x-s(t))\sqrt{\pi}e^{\alpha^{2}t}e^{-\alpha x}=|u_{0}|_{\alpha}\sup_{x>0}%
\exp(-\gamma t+\alpha^{2}t-\alpha x)\exp(\alpha\lbrack x-s(t)])\\
\leq|u_{0}|_{\alpha}\exp(\alpha^{2}-\gamma+\alpha V^{0})t
\end{gather*}
where we have used the bound $-s(t)\leq V^{0}t.$

To estimate the integral
\[
e^{\alpha^{2}t}e^{-\alpha x}\int_{\frac{(2t\alpha-x)}{2\sqrt{t}}}^{\infty}%
\exp(-\eta^{2}\}d\eta
\]
for $x<0,$ we apply Lemma \ref{errfun}. If the lower limit is larger than
unity then
\[
e^{\alpha^{2}t}e^{-\alpha x}\int_{\frac{(2t\alpha-x)}{2\sqrt{t}}}^{\infty}%
\exp(-\eta^{2}\}d\eta\leq\frac{1}{2}\exp(\alpha^{2}t-\alpha x-\frac
{(2t\alpha-x)^{2}}{4t})=\frac{1}{2}\exp(-\frac{x^{2}}{4t})
\]
If the lower limit is no larger than unity then
\begin{gather*}
e^{\alpha^{2}t}e^{-\alpha x}\int_{\frac{(2t\alpha-x)}{2\sqrt{t}}}^{\infty}%
\exp(-\eta^{2}\}d\eta\leq\frac{\sqrt{\pi}}{2}\exp(\alpha^{2}t-\alpha x)\\
=\frac{\sqrt{\pi}}{2}\exp(\frac{(2t\alpha-x)^{2}}{4t})\exp(-\frac{x^{2}}%
{4t})\leq\frac{e\sqrt{\pi}}{2}\exp(-\frac{x^{2}}{4t})
\end{gather*}
Thus, for $x<0$, we get
\begin{gather*}
\frac{e\sqrt{\pi}}{2}\frac{e^{-\gamma t}}{\sqrt{\pi}}|u_{0}|_{\alpha}%
\sup_{x<0}\exp(\alpha|x-s(t)|-\frac{x^{2}}{4t})\leq\\
\frac{e}{2}e^{-\gamma t}|u_{0}|_{\alpha}\max[\sup_{x<s(t)}\exp\{-\alpha
x+\alpha s(t)-\frac{x^{2}}{4t}\},\sup_{s(t)<x<0}\exp\{\alpha x-\alpha
s(t)-\frac{x^{2}}{4t}\}]\\
=\frac{e}{2}e^{-\gamma t}|u_{0}|_{\alpha}\max[\exp(\alpha^{2}t-\alpha
v_{0}t),\exp(V^{0}t\alpha)]
\end{gather*}
In the above estimate we used the elementary inequality:
\begin{align*}
-\alpha x+\alpha s(t)-\frac{x^{2}}{4t}  & =-(\frac{x}{\sqrt{4t}}+\sqrt{\alpha
t})^{2}+\alpha^{2}t+\alpha s(t)\leq\\
\alpha^{2}t+\alpha s(t)  & \leq\alpha^{2}t-\alpha v_{0}t
\end{align*}
(note that $s(t)\leq-v_{0}t$).

Collecting the estimates for all the cases ($x<0$, and $x>0$)
\begin{align*}
|u_{2}(.,t)|_{\alpha}  & =|u_{0}|_{\alpha}e^{-\gamma t}\max[\frac{e}{2}%
\exp(\alpha^{2}t-\alpha v_{0}t),\frac{e}{2}\exp(V^{0}t\alpha),\exp(\alpha
^{2}+\alpha V^{0})t]\\
& \leq\frac{e}{2}\exp[(-\gamma+\alpha^{2}+\alpha V^{0})t]\,|u_{0}|_{\alpha
}<2\exp[(-\gamma+\alpha^{2}+\alpha V^{0})t]\,|u_{0}|_{\alpha}%
\end{align*}
Thus, for any $\alpha$ we have obtained an a priori estimate on the
contribution from the initial data valid for all time. If $\alpha$ is
sufficiently small, $\alpha^{2}+\alpha V^{0}-\gamma<0,$ then the norm of the
contribution is exponentially decaying. We also note that for $\alpha
\rightarrow0$ the estimate has a limit and takes the form
\begin{equation}
|u_{2}(.,t)|_{0}\leq\frac{e}{2}\exp(-\gamma t)\,|u_{0}|_{0}\label{T2-solution}%
\end{equation}
\end{proof}

\subsection{Estimates for the solution: Contribution from the free interface}

\begin{proposition}
The $C_{\alpha}$-norm of the contribution from the free interface is uniformly
bounded for all time:
\[
|(T_{1}u_{0})(.,t)|_{\alpha}\leq V^{0}/\sqrt{\gamma},
\]
provided $\alpha<\alpha_{space}:=\min(v_{0}/4,\gamma/(2V^{0})).$
\end{proposition}

\begin{proof}
To estimate the free-interface contribution to the solution $T_{1}(t)u_{0}$
\textit{behind the interface} $x>s(t)$ we split the interval of integration
into two subsets: $\chi_{1}=\{\tau\in\lbrack0,t]:s(\tau)<(s(t)+x)/2\}$ and its
complement $\chi_{2}=\{\tau\in\lbrack0,t]:s(\tau)>(s(t)+x)/2\}$.
\[
\left|  T_{1}(t)u_{0}\right|  \leq\int_{0}^{t}G(x,t,s(\tau),\tau
)e^{-\gamma(t-\tau)}\,|v(\tau)|d\tau=\int\limits_{\chi_{1}}+\int
\limits_{\chi_{2}}=I_{1}+I_{2},
\]
For the first integral\ we have
\begin{align*}
I_{1}  & =\int\limits_{\chi_{1}}\frac{\exp[-(x-s(\tau))^{2}\frac{1}{4(t-\tau
)}]}{2\sqrt{\pi(t-\tau)}}e^{-\gamma(t-\tau)}\,|v(\tau)|d\tau\\
& \leq\frac{V^{0}}{2\sqrt{\pi}}\int\limits_{\chi_{1}}(t-\tau)^{-1/2}%
\exp[-(x-s(t))^{2}\frac{1}{16(t-\tau)}]e^{-\gamma(t-\tau)}d\tau\\
& \leq\frac{V^{0}}{2\sqrt{\pi}}\int\limits_{\chi_{1}}(t-\tau)^{-1/2}%
\exp[-(x-s(t))\frac{v_{0}}{8}]e^{-\gamma(t-\tau)}d\tau\\
& \leq\frac{V^{0}}{2\sqrt{\pi}}\exp[-(x-s(t))\frac{v_{0}}{4}]\int
\limits_{0}^{(x-s(t))/(2v_{0})}\eta^{-1/2}e^{-\gamma\eta}d\eta\\
& =\frac{V^{0}}{2\sqrt{\gamma}}\operatorname{erf}\left(  \sqrt{\gamma
\frac{x-s(t)}{2v_{0}}}\right)  \exp[-(x-s(t))\frac{v_{0}}{4}]\leq\frac{V^{0}%
}{2\sqrt{\gamma}}\exp[-(x-s(t))\frac{v_{0}}{4}]
\end{align*}
The following inequalities
\[
(x-s(\tau))^{2}\frac{1}{(t-\tau)}\leq(\frac{x-s(t)}{2})^{2}\frac{1}{(t-\tau
)}\leq(\frac{x-s(t)}{2})^{2}\frac{2v_{0}}{x-s(t)}%
\]
have been used to replace the exponent in the Gaussian kernel, which after
that gave the exponential decay factor. We note that the estimate has a
regular behavior at the limit $\gamma\rightarrow0$ giving the bound
\[
\frac{V^{0}\sqrt{x-s(t)}}{\sqrt{2v_{0}\pi}}\exp[-(x-s(t))\frac{v_{0}}{4}].
\]
It is a manifestation of the fact that the heat loss is immaterial in the
vicinity of the interface (cf. the next estimate which shows that the presence
of heat losses is essential for decay at large distances from the interface).

For the integral $I_{2}$ we use Lemma \ref{errfun} to obtain
\begin{gather}
I_{2}=%
%TCIMACRO{\dint \limits_{\chi_{2}}}%
%BeginExpansion
{\displaystyle\int\limits_{\chi_{2}}}
%EndExpansion
\frac{e^{-\gamma(t-\tau)}}{2\sqrt{\pi(t-\tau)}}\exp[-\dfrac{(x-s(\tau))^{2}%
}{4(t-\tau)}]\,|v(\tau)|d\tau\nonumber\\
\leq V^{0}\int\limits_{(x-s(t))/(2V^{0})}^{\infty}\frac{1}{2\sqrt{\pi\eta}%
}e^{-\gamma\eta}d\eta=\frac{V^{0}}{\sqrt{\pi}}\int\limits_{\sqrt
{(x-s(t))/(2V^{0})}}^{\infty}\exp(-\gamma\xi^{2})d\xi\nonumber\\
\leq\left\{
\begin{array}
[c]{l}%
\dfrac{V^{0}}{2\sqrt{\gamma\pi}}\exp(-\gamma(x-s(t))/(2V^{0})),\text{ for
}\gamma(x-s(t))/(2V^{0})>1\\
\dfrac{V^{0}}{2\sqrt{\gamma}},\text{ for }0\leq\gamma(x-s(t))/(2V^{0})<1
\end{array}
\right. \label{i2-function}%
\end{gather}
Thus for $x>s(t)$ we obtain
\begin{equation}
|T_{1}(t)u_{0}(x)|\leq\left\{
\begin{array}
[c]{c}%
\dfrac{V^{0}}{2\sqrt{\gamma}}\exp[-(x-s(t))\frac{v_{0}}{4}]+\\
\frac{V^{0}}{2\sqrt{\gamma}}\exp(-\gamma(x-s(t))/(2V^{0})),\text{ for }%
\gamma(x-s(t))/(2V^{0})>1\\
\dfrac{V^{0}}{2\sqrt{\gamma}}\exp[-(x-s(t))\frac{v_{0}}{4}]+\\
\dfrac{V^{0}}{2\sqrt{\gamma}}\text{ for }0<\gamma(x-s(t))/(2V^{0})<1
\end{array}
\right. \label{beh1}%
\end{equation}
We note that upon multiplication by the weight $\exp(\alpha(x-s(t)))$ the
right hand sides of the estimate (\ref{beh1}) are decaying exponentials,
provided
\begin{equation}
\alpha<\alpha_{space}:=\min(v_{0}/4,\gamma/(2V^{0})),\label{aspace}%
\end{equation}
that attain their maximum $V^{0}/\sqrt{\gamma}$ at $x-s(t)=0.$ \ Thus
\[
|T_{1}(t)u_{0}(x)|\leq V^{0}/\sqrt{\gamma}\text{ for }x>s(t)
\]

Ahead of the interface $x<s(t)$ we have:
\begin{align}
& |\int_{0}^{t}e^{-\gamma(t-\tau)}\frac{e^{-(x-s(\tau))^{2}/4(t-\tau)}}%
{\sqrt{4\pi(t-\tau)}}v(\tau)d\tau|\nonumber\\
& \leq\int_{0}^{t}e^{-\gamma(t-\tau)}\exp\{-\frac{[(x-s(t))+(s(t)-s(\tau
))]^{2}}{4(t-\tau)}\}\frac{|v(\tau)|d\tau}{\sqrt{4\pi(t-\tau)}}\nonumber\\
& \leq\exp\{-\frac{v_{0}|x-s(t)|}{2}-\frac{|x-s(t)|^{2}}{4t}\}\frac{V^{0}%
}{\sqrt{\pi}}\int_{0}^{t}\exp\{-\frac{(s(t)-s(\tau))^{2}}{4(t-\tau)}%
-\gamma(t-\tau)\}\frac{d\tau}{2\sqrt{(t-\tau)}}\nonumber\\
& \leq\exp\{-\frac{v_{0}|x-s(t)|}{2}-\frac{|x-s(t)|^{2}}{4t}\}\frac{V^{0}%
}{\sqrt{\pi}}\int_{0}^{t}\exp\{(-\frac{v_{0}^{2}}{4}-\gamma)(t-\tau
)\}\frac{d\tau}{2\sqrt{(t-\tau)}}\nonumber\\
& \leq\frac{V^{0}}{\sqrt{v_{0}^{2}+4\gamma}}\exp\{-\frac{v_{0}|x-s(t)|}%
{2}-\frac{|x-s(t)|^{2}}{4t}\}
\end{align}
Thus,
\[
|T_{1}(t)u_{0}(x)|\leq\frac{V^{0}}{\sqrt{v_{0}^{2}+4\gamma}}\exp\{-\frac
{v_{0}|x-s(t)|}{2}\}\text{ for }x<s(t)
\]
Now it is easy to obtain the estimate for the norm:
\begin{equation}
|(T_{1}u_{0})(.,t)|_{\alpha}\leq\max(\sup_{x<s(t)}[\frac{V^{0}}{\sqrt
{v_{0}^{2}+4\gamma}}\exp(-(\frac{v_{0}}{2}-\alpha)(s(t)-x))],\dfrac{V^{0}%
}{\sqrt{\gamma}})\leq\dfrac{V^{0}}{\sqrt{\gamma}}\label{T1-solution}%
\end{equation}
\end{proof}

By combining the estimates for the initial data contribution and that from the
free interface we arrive at the following result:

\begin{theorem}
\label{absorb}For $0\leq\alpha\leq\alpha_{space}$%
\[
|(Tu_{0})(.,t)|_{\alpha}\leq\dfrac{V^{0}}{\sqrt{\gamma}}+2\exp[(-\gamma
+\alpha^{2}+\alpha V^{0})t]\,|u_{0}|_{\alpha}%
\]
If, in addition,
\begin{equation}
\alpha<\alpha_{time}:=\dfrac{V^{0}}{2}(\sqrt{1+4\gamma/(V^{0})^{2}%
}-1)\label{atime}%
\end{equation}
(here $\alpha_{time}$ is the positive root of $\alpha^{2}+\alpha V^{0}%
-\gamma=0$) then the contribution from the initial data decays exponentially.
\end{theorem}

For the future use we combine the bounds:
\begin{equation}
\alpha_{\min}=\min(\alpha_{time},\alpha_{space})
\end{equation}

\begin{remark}
We note that for realistic problems $\gamma/V^{0}\ll1$ then $\alpha
_{time}\approx\gamma/V^{0},$ thus both bounds are of the same order and in
this case $\alpha_{\min}=\alpha_{space}.$
\end{remark}

\subsection{Estimates for the derivative: Contribution from initial data}

We also need estimates for the spatial derivative of the solution. We start
with the contribution from initial data.

\begin{proposition}
For $\alpha<\alpha_{time}$ the $C_{\alpha}$-norm of the derivative of the
contribution from the initial data decays exponentially in time:
\[
|(T_{2}u)_{x}(.,t)|_{\alpha}\leq|u_{0}|_{\alpha}(\frac{2}{\sqrt{t\pi}}%
+\frac{\alpha}{2})\exp[(\alpha^{2}+\alpha V^{0}-\gamma)t].
\]
\end{proposition}

\begin{proof}
First we split the integral
\begin{align*}
& |(T_{2}u)_{x}(.,t)|_{a}=\sup\left\{  \omega_{\alpha}(x-s(t))\;|\int
_{-\infty}^{\infty}e^{-\gamma t}G_{x}(x,t,\xi,0)u_{0}(\xi)d\xi|\right\} \\
& =\sup\left\{  e^{-\gamma t}\omega_{\alpha}(x-s(t))|\int_{-\infty}^{\infty
}\frac{1}{\omega_{\alpha}(\xi)}G_{x}(x,t,\xi,0)u_{0}(\xi)\omega_{\alpha}%
(\xi)d\xi|\right\} \\
& =\frac{e^{-\gamma t}}{2\sqrt{t\pi}}|u_{0}|_{\alpha}\,\sup\{\omega_{\alpha
}(x-s(t))[\int_{0}^{\infty}e^{-\alpha\xi}\left|  \frac{\partial}{\partial
x}\exp(-\frac{(x-\xi)^{2}}{4t})\right|  d\xi\\
& +\int_{-\infty}^{0}e^{\alpha\xi}\left|  \frac{\partial}{\partial x}%
\exp(-\frac{(x-\xi)^{2}}{4t})\right|  d\xi]\}
\end{align*}
Each of the integrals should be estimated separately for $x<0$ and $x>0$ and
the maximum of the estimates should be chosen for the estimate of the norm$.$
At the same time it is easy to see that the two integrals can be transformed
to each other through the change $x\rightarrow-x,$ therefore it suffices to
estimate only one of them and double the result.

We proceed as follows. For $x>0$ we integrate by parts
\begin{gather*}
\int_{0}^{\infty}e^{-\alpha\xi}\left|  \frac{\partial}{\partial x}\exp
(-\frac{(x-\xi)^{2}}{4t})\right|  d\xi=\int_{0}^{\infty}e^{-\alpha\xi}\left|
\frac{\partial}{\partial\xi}\exp(-\frac{(x-\xi)^{2}}{4t})\right|  d\xi\\
=\int_{0}^{x}e^{-\alpha\xi}\frac{\partial}{\partial\xi}\exp(-\frac{(x-\xi
)^{2}}{4t})d\xi-\int_{x}^{\infty}e^{-\alpha\xi}\frac{\partial}{\partial\xi
}\exp(-\frac{(x-\xi)^{2}}{4t})d\xi\\
=\alpha(\int_{0}^{x}e^{-\alpha\xi}\exp(-\frac{(x-\xi)^{2}}{4t})d\xi-\int
_{x}^{\infty}e^{-\alpha\xi}\exp(-\frac{(x-\xi)^{2}}{4t})d\xi)+2e^{-\alpha
x}-\exp(-\frac{x^{2}}{4t})\\
=\alpha\left\{  \int_{-x}^{0}e^{-\alpha(x+\eta)}\exp(-\frac{\eta^{2}}%
{4t})d\eta-\int_{0}^{\infty}e^{-\alpha(x+\eta)}\exp(-\frac{\eta^{2}}{4t}%
)d\eta\right\}  +2e^{-\alpha x}-\exp(-\frac{x^{2}}{4t})\\
\leq\alpha\exp(\alpha^{2}t-\alpha x)\int_{-x}^{0}\exp(-\frac{\eta^{2}%
+4\alpha\eta t+4\alpha^{2}t^{2}}{4t})d\eta+2e^{-\alpha x}\\
\leq\alpha\exp(\alpha^{2}t-\alpha x)\int_{-\infty}^{\infty}\exp(-\frac
{(\eta+2\alpha t)^{2}}{4t})d\eta+2e^{-\alpha x}\\
=\alpha\exp(\alpha^{2}t-\alpha x)\sqrt{t\pi}+2e^{-\alpha x}%
\end{gather*}

For $x<0$ we integrate by parts to obtain
\begin{gather*}
\int_{0}^{\infty}e^{-\alpha\xi}\left|  \frac{\partial}{\partial x}\exp
(-\frac{(x-\xi)^{2}}{4t})\right|  d\xi=\int_{0}^{\infty}e^{-\alpha\xi}\left|
\frac{\partial}{\partial\xi}\exp(-\frac{(x-\xi)^{2}}{4t})\right|  d\xi=\\
-\int_{0}^{\infty}e^{-\alpha\xi}\frac{\partial}{\partial\xi}\exp(-\frac
{(x-\xi)^{2}}{4t})d\xi=-\alpha\int_{-x}^{\infty}e^{-\alpha(x+\xi)}\exp
(-\frac{\xi^{2}}{4t})d\xi+\exp(-\frac{x^{2}}{4t})\leq\exp(-\frac{x^{2}}{4t})
\end{gather*}

Now we are ready to estimate the norm:
\begin{align*}
|(T_{2}u)_{x}(.,t)|_{\alpha}  & \leq\frac{e^{-\gamma t}}{2\sqrt{t\pi}}%
|u_{0}|_{\alpha}\max[\sup_{x>0}\left\{  \omega_{\alpha}(x-s(t))(\alpha
\sqrt{t\pi}\exp(\alpha^{2}t-\alpha x)+2e^{-\alpha x})\right\}  ,\\
& \sup_{x<0}\left\{  \omega_{\alpha}(x-s(t))\exp(-\frac{x^{2}}{4t}))\right\}
]
\end{align*}
For the term with $x>0$ in the above estimate, we have
\begin{align*}
& \frac{e^{-\gamma t}}{2\sqrt{t\pi}}|u_{0}|_{\alpha}\sup_{x>0}\left\{
\omega_{\alpha}(x-s(t))[\alpha\sqrt{t\pi}\exp(\alpha^{2}t-\alpha
x)+2e^{-\alpha x})]\right\} \\
& \leq\frac{e^{-\gamma t}}{2\sqrt{t\pi}}|u_{0}|_{\alpha}\exp(-\alpha
s(t))[2+\alpha\sqrt{t\pi}\exp(\alpha^{2}t)]\\
& \leq|u_{0}|_{\alpha}\left\{  \frac{1}{\sqrt{t\pi}}\exp[(\alpha V^{0}%
-\gamma)t)]+\frac{\alpha}{2}\exp[(\alpha^{2}+\alpha V^{0}-\gamma)t]\right\}
\end{align*}

For $x<0,$ we have
\begin{gather*}
|(T_{2}u)_{x}(.,t)|_{\alpha}\leq\frac{e^{-\gamma t}}{2\sqrt{t\pi}}%
|u_{0}|_{\alpha}\sup_{x<0}\left\{  \omega_{\alpha}(x-s(t))\exp(-\frac{x^{2}%
}{4t}))\right\}  ]\\
\leq\frac{e^{-\gamma t}}{2\sqrt{t\pi}}|u_{0}|_{\alpha}\max[\sup_{x<s(t)}%
\exp\{-\alpha x+\alpha s(t)-\frac{x^{2}}{4t}\},\sup_{s(t)<x<0}\exp\{\alpha
x-\alpha s(t)-\frac{x^{2}}{4t}\}]\\
\leq\frac{1}{2\sqrt{t\pi}}|u_{0}|_{\alpha}\max\left\{  \exp[(\alpha^{2}-\alpha
v_{0}-\gamma)t],\exp[(\alpha V^{0}-\gamma)t]\right\}
\end{gather*}
In the above estimate we have replaced $\alpha x-\dfrac{x^{2}}{4t}$ by its
maximum $\alpha^{2}t.$ We collect the estimate for $x>0$ and $x<0$ to obtain
\begin{align}
|(T_{2}u)_{x}(.,t)|_{\alpha}  & \leq|u_{0}|_{\alpha}\left\{  \frac{1}%
{\sqrt{t\pi}}\exp[(\alpha V^{0}-\gamma)t)]+(\frac{1}{\sqrt{t\pi}}+\frac
{\alpha}{2})\exp[(\alpha^{2}+\alpha V^{0}-\gamma)t]\right\} \nonumber\\
& \leq|u_{0}|_{\alpha}(\frac{2}{\sqrt{t\pi}}+\frac{\alpha}{2})\exp[(\alpha
^{2}+\alpha V^{0}-\gamma)t].\label{T2-derivative}%
\end{align}
Thus, for any $\alpha$ we have obtained an a priori estimate on the
contribution from the initial data valid for all time. If $0\leq\alpha
<\alpha_{time}$ then the norm of the contribution is exponentially decaying.
\end{proof}

\subsection{Estimates for the derivative: Contribution from the interface}

\begin{proposition}
The $C_{\alpha}$-norm of the derivative of the contribution from the free
interface is uniformly bounded for all time:
\[
|(T_{1}u)_{x}(.,t)|_{\alpha}\leq\mathcal{M}(v_{0},V^{0},\alpha,\gamma)
\]
provided that
\[
0\leq\alpha<\min(\frac{v_{0}}{8},\frac{\gamma}{2V^{0}}):=\alpha_{\min}%
^{\prime}%
\]
\end{proposition}

\begin{proof}
\textbf{(Ahead of the interface)}. The estimate ahead of the front $x\leq
s(t)$ is treated as follows. We consider separately two cases: $|s(t)-x|>1$
and $|s(t)-x|\leq1$.

For the case $|s(t)-x|>1$%
\begin{gather*}
|(T_{1}u)_{x}(x,t)|=|\int_{0}^{t}e^{-\gamma(t-\tau)}\frac{x-s(\tau)}%
{2(t-\tau)}\frac{e^{-(x-s(\tau))^{2}/4(t-\tau)}}{\sqrt{4\pi(t-\tau)}}%
v(\tau)d\tau|\\
\leq|\int_{0}^{t}\frac{(x-s(\tau))^{2}}{2(t-\tau)(x-s(\tau))}e^{-(x-s(\tau
))^{2}/8(t-\tau)}\times e^{-(x-s(\tau))^{2}/8(t-\tau)}\frac{v(\tau)d\tau
}{\sqrt{4\pi(t-\tau)}}|\\
\leq|\int_{0}^{t}\frac{4/e}{s(t)-s(\tau)}\exp\{-\frac{[(x-s(t))+(s(t)-s(\tau
))]^{2}}{8(t-\tau)}\}\frac{v(\tau)d\tau}{\sqrt{\pi(t-\tau)}}|\\
\leq\frac{4V^{0}}{v_{0}e\sqrt{\pi}}\int_{0}^{t}(t-\tau)^{-3/2}\exp
[-\frac{(x-s(t))^{2}}{8(t-\tau)}-\frac{v_{0}}{4}|x-s(t)|-\frac{v_{0}^{2}}%
{8}(t-\tau)]d\tau\\
\leq\frac{4V^{0}e^{-v_{0}|x-s(t)|/4}}{e\sqrt{\pi}v_{0}|s(t)-x|}\int_{0}%
^{t}|s(t)-x|(t-\tau)^{-3/2}\exp[-\frac{(x-s(t))^{2}}{8(t-\tau)}]d\tau\\
\leq\frac{\sqrt{256}V^{0}e^{-v_{0}|x-s(t)|/4}}{e\sqrt{\pi}v_{0}|s(t)-x|}%
\int_{0}^{\infty}e^{-\eta^{2}}d\eta\leq\frac{8}{e}\frac{V^{0}e^{-v_{0}%
|x-s(t)|/4}}{v_{0}|s(t)-x|}%
\end{gather*}
In the last estimate we used the following simple observations: $\xi e^{-\xi
}\leq1/e,$ for $\xi=\dfrac{(x-s(\tau))^{2}}{4(t-\tau)}>0,$ $|s(\tau
)-x|>|s(t)-x|$, \ $|s(\tau)-x|>|s(t)-s(\tau)|>v_{0}|t-\tau|$ and substitution
$\eta=|s(t)-x|(t-\tau)^{-1/2}$ to obtain the error function integral.

For the less involved case $|s(t)-x|\leq1$ we proceed as follows
\begin{gather*}
|(T_{1}u)_{x}(x,t)|\leq|\int_{0}^{t}\frac{x-s(\tau)}{2(t-\tau)}\frac
{e^{-(x-s(\tau))^{2}/4(t-\tau)}}{\sqrt{4\pi(t-\tau)}}v(\tau)d\tau|\\
\leq|\int_{0}^{t}\frac{|x-s(t)|+|s(t)-s(\tau)|}{2(t-\tau)}\frac{e^{-(x-s(\tau
))^{2}/4(t-\tau)}}{\sqrt{4\pi(t-\tau)}}v(\tau)d\tau|\\
\leq\frac{V^{0}}{\sqrt{\pi}}\int_{0}^{t}\frac{|s(t)-x|(t-\tau)^{-3/2}}%
{4}e^{-(x-s(t))^{2}/4(t-\tau)}d\tau+\frac{(V^{0})^{2}}{4\sqrt{\pi}}\int
_{0}^{t}\frac{e^{-(s(t)-s(\tau))^{2}/4(t-\tau)}}{\sqrt{(t-\tau)}}d\tau\\
\leq\frac{V^{0}}{\sqrt{\pi}}\int_{0}^{\infty}e^{-\eta^{2}}d\eta+\frac
{(V^{0})^{2}}{4\sqrt{\pi}}\int_{0}^{t}\frac{e^{-v_{0}^{2}(t-\tau)/4}}%
{\sqrt{(t-\tau)}}d\tau\\
\leq\frac{V^{0}}{2}+\frac{(V^{0})^{2}}{2\sqrt{\pi}}\int_{0}^{\infty}%
e^{-v_{0}^{2}\eta^{2}/4}d\eta\leq\frac{V^{0}}{2}(1+\frac{V^{0}}{v_{0}})
\end{gather*}
Thus for $x<s(t)$ \ we obtain
\begin{equation}
|(T_{1}u)_{x}(x,t)|\leq\left\{
\begin{array}
[c]{c}%
\dfrac{8V^{0}\exp(-v_{0}|x-s(t)|/4)}{ev_{0}|s(t)-x|}\text{ for }s(t)-x>1\\
\dfrac{V^{0}}{2}(1+\dfrac{V^{0}}{v_{0}})\text{ for }0<s(t)-x<1
\end{array}
\right. \label{nehi}%
\end{equation}
\end{proof}

\begin{remark}
The proof above shows that if $|s(t)-s(\tau)|\geq v_{0}|t-\tau|$ which holds
if the basic assumption on the kinetics in (\ref{kinetics}) is satisfied, then
the derivative ahead of the interface $x<s(t)$ decays exponentially
\begin{equation}
|(T_{1}u)_{x}(x,t)|\leq\frac{Ce^{-v_{0}|x-s(t)|/4}}{|s(t)-x|}\left\|
v\right\|  _{C[0,t]}\label{expdec}%
\end{equation}
The exponent $-v_{0}/4$ can be improved to $-v_{0}/(2+\varepsilon)$ (at the
price of increasing $C$).
\end{remark}

\begin{proof}
\textbf{(Behind the interface)}. Now we consider the domain behind the
interface, $x>s(t)$. We split the interval of integration into two subsets:
$\chi_{1}=\{\tau\in\lbrack0,t]:s(\tau)<(s(t)+x)/2\}$ and its compliment
$\chi_{2}=\{\tau\in\lbrack0,t]:s(\tau)>(s(t)+x)/2\}$.
\[
|(T_{1}u)_{x}(x,t)|\leq\int_{0}^{t}\dfrac{|x-s(\tau)|}{2(t-\tau)}%
G(x,t,s(\tau),\tau)e^{-\gamma(t-\tau)}\,|v(\tau)|d\tau=\int\limits_{\chi_{1}%
}+\int\limits_{\chi_{2}}=I_{1}+I_{2},
\]

For the first integral\ we have
\begin{gather*}
I_{1}=\frac{1}{4\sqrt{\pi}}\int\limits_{\chi_{1}}\dfrac{|x-s(\tau)|}%
{(t-\tau)^{3/2}}\exp[-\frac{(x-s(\tau))^{2}}{4(t-\tau)}]e^{-\gamma(t-\tau
)}|v(\tau)|d\tau|\\
\leq\frac{V^{0}}{4\sqrt{\pi}}\int\limits_{\chi_{1}}\dfrac{|x-s(t)|}%
{(t-\tau)^{3/2}}\exp[-\frac{(x-s(t))^{2}}{16(t-\tau)}]e^{-\gamma(t-\tau)}%
d\tau\\
\leq\frac{V^{0}}{4\sqrt{\pi}}\int\limits_{0}^{(x-s(t))/(2v_{0})}%
\dfrac{|x-s(t)|}{(t-\tau)^{3/2}}\exp[-\frac{(x-s(t))^{2}}{16(t-\tau
)}]e^{-\gamma(t-\tau)}d(t-\tau)=\frac{V^{0}}{2\sqrt{\pi}}\int\limits_{\sqrt
{(x-s(t))2v_{0}}}^{\infty}e^{-\eta^{2}/16}d\eta
\end{gather*}
where $\eta=(x-s(t))/\sqrt{t-\tau}.$

To estimate the last integral we apply Lemma \ref{errfun} to obtain
\[
I_{1}\leq\left\{
\begin{array}
[c]{l}%
\dfrac{2V^{0}}{\sqrt{\pi}}\exp(-(x-s(t))v_{0}/8),\text{ for }(x-s(t))2v_{0}%
>1\\
2V^{0},\text{ for }0\leq(x-s(t))2v_{0}<1
\end{array}
\right.
\]
For the estimate on the $C_{\alpha}$-norm, the function will be multiplied by
the weight. We note that upon multiplication by $\exp(\alpha(x-s(t)))$ the
right hand sides of the top estimate is a decaying exponential, provided
$\alpha<v_{0}/8,$ that attain its maximum $2V^{0}/\sqrt{\pi}$ at $x-s(t)=0, $
while the bottom term is bounded by $2V^{0}\exp(\alpha/(v_{0}/8))<2eV^{0}.$
Therefore the contribution of $I_{1}$ into the norm is bounded above, for
example, by $6V^{0}$: $I_{1}<6V^{0}.$

For the integral $I_{2}$%
\[
I_{2}=\int\limits_{\chi_{2}}\dfrac{|x-s(\tau)|}{2(t-\tau)}\frac{\exp
[-\dfrac{(x-s(\tau))^{2}}{4(t-\tau)}]}{2\sqrt{\pi(t-\tau)}}e^{-\gamma(t-\tau
)}\,|v(\tau)|d\tau
\]
we use simple geometric considerations that show that in the domain $\chi_{2}$
\ if $s(\tau)>x$ then
\[
\dfrac{|x-s(\tau)|}{2(t-\tau)}\leq\dfrac{|s(t)-s(\tau)|}{2(t-\tau)}\leq
\frac{V^{0}}{2}%
\]
while for $s(\tau)\leq x$%
\[
\dfrac{|x-s(\tau)|}{2(t-\tau)}\leq\frac{1}{2}|x-\frac{s(t)+x}{2}%
|/(|s(t)-\frac{s(t)+x}{2}|/V^{0})=\frac{V^{0}}{2}.
\]
Therefore the estimate for the derivative in this case reduces to the estimate
for the function itself (\ref{i2-function}) and yields:
\[
I_{2}\leq\frac{V^{0}}{2}\left\{
\begin{array}
[c]{l}%
\frac{V^{0}}{2\sqrt{\gamma\pi}}\exp(-\gamma(x-s(t))/(2V^{0})),\text{ for
}\gamma(x-s(t))/(2V^{0})>1\\
\frac{V^{0}}{2\sqrt{\gamma}},\text{ for }0\leq\gamma(x-s(t))/(2V^{0})<1
\end{array}
\right.
\]
Similarly to the argument for $I_{1},$ one can see that if $\alpha
<\gamma/(2V^{0})$ then the contribution of $I_{2}$ into the $C_{\alpha}$-norm
is bounded by $V^{0}e/(2\sqrt{\gamma})$.

Thus for $x>s(t)$ \ we obtain
\begin{equation}
\sup_{x>s(t)}|(T_{1}u)_{x}(x,t)\exp(\alpha(x-s(t))|\leq V^{0}e/(2\sqrt{\gamma
})+6V^{0}<V^{0}(2/\sqrt{\gamma}+6)\label{behi}%
\end{equation}
if \ $\alpha<\min(\frac{v_{0}}{8},\frac{\gamma}{2V^{0}}).$

Similarly, if $\alpha<v_{0}/4$ then by employing (\ref{nehi}) we see that
\[
\sup_{x<s(t)}|(T_{1}u)_{x}(x,t)\exp(-\alpha(x-s(t))|\leq\max[\dfrac{8V^{0}%
}{ev_{0}},\dfrac{V^{0}}{2}(1+\dfrac{V^{0}}{v_{0}})e^{\alpha}]
\]
Finally, for the norm we get
\begin{equation}
|(T_{1}u)_{x}(.,t)|_{\alpha}\leq V^{0}\max[\dfrac{8}{ev_{0}},\dfrac{1}%
{2}(1+\dfrac{V^{0}}{v_{0}})e^{\alpha},2/\sqrt{\gamma}+6]:=\mathcal{M}%
(v_{0},V^{0},\alpha,\gamma)\label{T1-derivative}%
\end{equation}
The estimate holds if
\begin{equation}
0\leq\alpha<\min(\frac{v_{0}}{8},\frac{\gamma}{2V^{0}}):=\alpha_{\min}%
^{\prime}\label{amin}%
\end{equation}
(Recall that $T_{1}u$ is the contribution from the free interface and
therefore the absolute bound on its derivative is independent of the initial data).
\end{proof}

\begin{remark}
\label{4+epsilon} The choice of the factor $1/2$ for the split of the domain
of integration into two parts above is rather arbitrary; by choosing a factor
approaching $1$ we can improve the exponent of decay in the subsequent result,
the corresponding coefficient though will increase. Consequently, the value
$\dfrac{v_{0}}{8}$ in the definition of $\alpha_{\min}^{\prime}$ can be
improved to become $\dfrac{v_{0}}{4+\varepsilon}. $
\end{remark}

\subsection{A priori estimates in $C_{\alpha}^{1}$ and $H_{\alpha}^{1}$}

We collect the estimates (\ref{T2-derivative})-(\ref{T1-derivative}) to obtain
the following result:

\begin{theorem}
\label{c1alpha}
\[
\left|  (Tu_{0})_{x}(.,t)\right|  _{\alpha}\leq\mathcal{M}+|u_{0}|_{\alpha
}(\frac{2}{\sqrt{t\pi}}+\frac{\alpha}{2})\exp[(\alpha V^{0}+\alpha^{2}%
-\gamma)t]
\]
\end{theorem}

\begin{remark}
It is easy to demonstrate via integration by parts that if the initial data
$u_{0}\in C_{\alpha}^{1}$ and satisfy the compatibility condition
$[(u_{0})_{x}]_{x=0}=g(u_{0}(0))\;$ then the $1/\sqrt{t}$ singularity in the
above estimate will not take place. The estimate in this case reduces to the
estimate for the solution through the derivative of the initial conditions.
\end{remark}

We note that a weaker result holds for Hilbert norms.

\begin{theorem}
Let $u_{0}\in C_{\alpha}$ then the solution satisfies $u(.,t)\in
C((0,\infty),H_{\beta}^{1})$ where $\beta<\alpha$ and
\begin{align*}
\left\|  (Tu_{0})(.,t)\right\|  _{\beta,1}  & \leq\frac{1}{\sqrt{\alpha-\beta
}}\{\mathcal{M}+|u_{0}|_{\alpha}(\frac{2}{\sqrt{t\pi}}+\frac{\alpha}{2}%
)\exp[(\alpha V^{0}+\alpha^{2}-\gamma)t]\\
& +\dfrac{V^{0}}{\sqrt{\gamma}}+2\,|u_{0}|_{\alpha}\exp[(-\gamma+\alpha
^{2}+\alpha V^{0})t]\}
\end{align*}
\end{theorem}

\begin{proof}
The proof is very simple since both $Tu_{0}(.,t)$ and $(Tu_{0})_{x}(.,t)\in
H_{\beta}$ \ by virtue of the imbedding estimate (\ref{beta-alpha}).
\end{proof}

\section{Well-posedness}

In this section we prove that solutions of the free boundary problem
(\ref{he})-(\ref{jc}) depend continuously on initial data. This result is used
in the sequel to demonstrate smoothness of the elements of the attractor.

\begin{theorem}
\label{well-pose}In $C_{\alpha},$ $0\leq\alpha<\alpha_{space},$ the solutions
of the problem depend on initial conditions uniformly continuously. More
precisely, if $\{u(x,t),s(t)\},\ \{\tilde{u}(x,t),\tilde{s}(t)\},\ 0<t<\sigma$
are solutions with initial data $u^{0},\tilde{u}^{0}\in C_{\alpha}$, where
$\sigma>0$ depends only on the norm $|u^{0}|_{\alpha}$ and $|\tilde{u}%
^{0}|_{\alpha},$ then for $0<t<\sigma$
\begin{equation}
\sup_{0<t<\sigma}|V(t)-\tilde{V}(t)|<c|u^{0}-\tilde{u}^{0}|_{\alpha
},\label{v-estim}%
\end{equation}%
\begin{equation}
|u(.-s(t),t)-\tilde{u}(.-\tilde{s}(t),t)|_{\alpha}<c|u^{0}-\tilde{u}%
^{0}|_{\alpha}\label{u-estim}%
\end{equation}
\end{theorem}

\begin{remark}
We state and prove continuous dependence on initial conditions only locally in
time. The argument extending this result to any fixed time is based on the a
priori estimates and follows closely the proof of global existence.
\end{remark}

\begin{remark}
For simplicity of presentation we include the proof only for the uniform norm
$\alpha=0.$ The modifications for the case $\alpha>0$ are rather routine but
somewhat lengthy and follow along the similar lines. Everywhere in the proof
below we use the notation $\left\|  .\right\|  =|.|_{0}$
\end{remark}

\begin{proof}
The proof consists of two parts. First we establish continuity of the
interface velocity. We will establish first the estimate in (\ref{v-estim})
and then use it to derive (\ref{u-estim}). Let $v$ and $\tilde{v}$ be
solutions of the integral equation in (\ref{ie}) with initial data $u^{0}$ and
$\tilde{u}^{0}$ respectively.
\begin{align*}
\left\|  v-\tilde{v}\right\|   & =\left\|  g\left(  e^{-\gamma t}\int
_{-\infty}^{\infty}G(s(t),t,\xi,0)u_{0}(\xi)d\xi-\int_{0}^{t}G(s(t),t,s(\tau
),\tau)e^{-\gamma(t-\tau)}v(\tau)d\tau\right)  \right. \\
& \left.  -g\left(  e^{-\gamma t}\int_{-\infty}^{\infty}G(\tilde{s}%
(t),t,\xi,0)\tilde{u}_{0}(\xi)d\xi-\int_{0}^{t}G(\tilde{s}(t),t,\tilde{s}%
(\tau),\tau)e^{-\gamma(t-\tau)}\tilde{v}(\tau)d\tau\right)  \right\| \\
& \leq Le^{-\gamma t}\left\|  \int_{-\infty}^{\infty}G(s(t),t,\xi,0)u_{0}%
(\xi)d\xi-\int_{-\infty}^{\infty}G(\tilde{s}(t),t,\xi,0)\tilde{u}_{0}(\xi
)d\xi\right\| \\
& +L\left\|  \int_{0}^{t}G(s(t),t,s(\tau),\tau)e^{-\gamma(t-\tau)}v(\tau
)d\tau-\int_{0}^{t}G(\tilde{s}(t),t,\tilde{s}(\tau),\tau)e^{-\gamma(t-\tau
)}\tilde{v}(\tau)d\tau\right\|
\end{align*}
To continue the estimate we employ a ''coordinate descent'':
\begin{align*}
& \leq Le^{-\gamma t}\left\|  \int_{-\infty}^{\infty}(G(s(t),t,\xi
,0)-G(\tilde{s}(t),t,\xi,0))u_{0}(\xi)d\xi\right\| \\
& +Le^{-\gamma t}\left\|  \int_{-\infty}^{\infty}G(\tilde{s}(t),t,\xi
,0)(u_{0}(\xi)-\tilde{u}_{0}(\xi))d\xi\right\| \\
& +L\left\|  \int_{0}^{t}G(s(t),t,s(\tau),\tau)e^{-\gamma(t-\tau)}%
(v(\tau)-\tilde{v}(\tau))d\tau\right\| \\
& +L\left\|  \int_{0}^{t}[G(s(t),t,s(\tau),\tau)-G(\tilde{s}(t),t,\tilde
{s}(\tau),\tau)]e^{-\gamma(t-\tau)}\tilde{v}(\tau)d\tau\right\| \\
& :=Le^{-\gamma t}D_{1}+Le^{-\gamma t}D_{2}+LD_{3}+LD_{4}%
\end{align*}

To estimate the first summand
\[
D_{1}=\int_{-\infty}^{\infty}(G(s(t),t,\xi,0)-G(\tilde{s}(t),t,\xi
,0))u_{0}(\xi)d\xi:=\int_{-\infty}^{\infty}\delta Gu_{0}d\xi
\]
we note that by the mean value theorem,
\[
\delta G=(s-\tilde{s})\frac{\partial G}{\partial x}(s^{\prime}-\xi
,t,0,0)=(s-\tilde{s})\frac{s^{\prime}-\xi}{2t}G(s^{\prime}-\xi,t,0,0)
\]
where
\[
s^{\prime}=s^{\prime}(t,\xi),\quad s(t)\leq s^{\prime}\leq\tilde{s}(t).
\]
Thus
\begin{align*}
|s^{\prime}-\xi|G(s^{\prime}-\xi,t,0,0)  & =(2\sqrt{\pi})^{-1}|s^{\prime}%
-\xi|t^{-1/2}e^{-(s^{\prime}-\xi)^{2}/4t}\\[0.08in]
& =(2\sqrt{\pi})^{-1}(8t)^{1/2}\frac{|s^{\prime}-\xi|}{(8t)^{1/2}%
}e^{-(s^{\prime}-\xi)^{2}/8t}2^{1/2}(2t)^{-1/2}e^{-(s^{\prime}-\xi)^{2}%
/8t}\\[0.08in]
& \leq4t^{1/2}c_{1}G(s^{\prime}-\xi,2t,0,0)\leq4c_{1}t^{1/2}G(s-\xi,2t,0,0)
\end{align*}
where $c_{1}=(2\sqrt{\pi})^{-1}\max(xe^{-x^{2}})$. Therefore
\begin{align*}
|D_{1}|  & =|\int_{-\infty}^{\infty}\delta Gu^{0}d\xi|\leq\frac{s-\tilde{s}%
}{2t}4c_{1}t^{1/2}\int_{-\infty}^{\infty}G(s-\xi,2t,0,0)|u^{0}(\xi
)|d\xi\\[0.08in]
& \leq2c_{1}t^{1/2}\frac{1}{t}\int_{0}^{t}|v-\tilde{v}|d\tau\int_{-\infty
}^{\infty}G(s-\xi,2t,0,0)|u^{0}(\xi)|d\xi\\[0.08in]
& \leq2c_{1}t^{1/2}\sup|u^{0}(\xi)|\,\Vert v-\tilde{v}\Vert.
\end{align*}

Obviously,
\[
|D_{2}|\leq\int_{-\infty}^{\infty}G(\tilde{s}(t),t,\xi,0)|u_{0}(\xi)-\tilde
{u}_{0}(\xi)|d\xi\leq\Vert u_{0}-\tilde{u}_{0}\Vert
\]
For the estimate of $D_{3}$ we replace the exponentials by $1$ and integrate
to obtain$:$
\[
|D_{3}|\leq\frac{1}{\sqrt{\pi}}t^{1/2}\Vert v-\tilde{v}\Vert
\]
Finally, for $D_{4}$%
\begin{align*}
& \left|  \int_{0}^{t}[G(s(t),t,s(\tau),\tau)-G(\tilde{s}(t),t,\tilde{s}%
(\tau),\tau)]e^{-\gamma(t-\tau)}\tilde{v}(\tau)d\tau\right| \\
& \leq\left\|  \tilde{v}\right\|  \int_{0}^{t}|\Delta G|e^{-\gamma(t-\tau
)}\tilde{v}(\tau)d\tau\tilde{v}(\tau)
\end{align*}
The estimations are quite elementary and are based on the mean value theorem.
First we note that
\begin{align*}
|\Delta G|  & :=|G(s(t),t,s(\tau),\tau)-G(\tilde{s}(t),t,\tilde{s}(\tau
),\tau)|\\[0.08in]
& =|G(s(t)-s(\tau),t-\tau,0,0)-G(\tilde{s}(t)-\tilde{s}(\tau),t-\tau
,0,0)|\\[0.08in]
& =|s(t)-s(\tau)-(\tilde{s}(t)-\tilde{s}(\tau))|\,|\frac{\partial G}{\partial
x}(s^{\prime},t-\tau,0,0)|\\[0.08in]
& =|\frac{s(t)-\tilde{s}(t)-(s(\tau)-\tilde{s}(\tau))}{2(t-\tau)}%
|\,|s^{\prime}G(s^{\prime},t-\tau,0,0)|\\[0.08in]
& =\frac{1}{2}|\frac{ds}{dt}(\tau^{\prime})-\frac{d\tilde{s}}{dt}(\tau
^{\prime})|\,|s^{\prime}G(s^{\prime},t-\tau,0,0)|\\
& \leq\frac{1}{2}\Vert v-\tilde{v}\Vert|s^{\prime}G(s^{\prime},t-\tau,0,0)|
\end{align*}
where $\tau\leq\tau^{\prime}\leq t$ and $s^{\prime}$ is between $\tilde
{s}(t)-\tilde{s}(\tau)$ and $s(t)-s(\tau)$. Since
\[
|s^{\prime}|\leq\max\left\{  |\tilde{s}(t)-\tilde{s}(\tau)|,|s(t)-s(\tau
)|\right\}  \leq V_{0}(t-\tau).
\]
and $|G|\leq C_{0}(t-\tau)^{-1/2}$ we get the estimate
\[
|\Delta G|\leq C_{0}V_{0}\Vert v-\tilde{v}\Vert(t-\tau)^{1/2}%
\]
Thus, for $D_{4}$ we get
\[
|D_{4}|\leq C_{4}\Vert v-\tilde{v}\Vert t^{3/2}%
\]
We collect the estimates for\ $D_{1}$-\ $D_{4}$\ to obtain
\[
\left\|  v-\tilde{v}\right\|  \leq\Vert u_{0}-\tilde{u}_{0}\Vert+C_{1}%
t^{1/2}\left\|  u^{0}\right\|  \,\Vert v-\tilde{v}\Vert+C_{4}\Vert v-\tilde
{v}\Vert t^{3/2}+C_{2}t^{1/2}\Vert v-\tilde{v}\Vert
\]
From this inequality we see that for $t<\sigma,$ where\ $\sigma$\ is small
enough
\[
\left\|  v-\tilde{v}\right\|  \leq\frac{\Vert u_{0}-\tilde{u}_{0}\Vert
}{1-C_{1}t^{1/2}\left\|  u^{0}\right\|  -C_{4}t^{3/2}-C_{2}t^{1/2}}\leq2\Vert
u_{0}-\tilde{u}_{0}\Vert
\]
We note that the value of\ $\sigma$\ depends only on $\left\|  u^{0}\right\|
$ therefore this result can be interpreted as uniformly continuous dependence
of $v$ on the initial data on any ball $\left\|  u^{0}\right\|  \leq R.$

Thus, the estimate (\ref{v-estim}) has been demonstrated. Next we note that
\begin{align}
& |u(x-s(t),t)-\tilde{u}(x-\tilde{s}(t),t)|\label{u-util}\\
& \leq|u(x-\tilde{s}(t),t)-\tilde{u}(x-\tilde{s}%
(t),t)|+|u(x-s(t),t)-u(x-\tilde{s}(t),t)|\nonumber
\end{align}
The second term in the inequality above can be estimated via the mean value
theorem:
\begin{align*}
|u(x-s(t),t)-u(x-\tilde{s}(t),t)|  & =|u_{x}(x-s^{\prime},t)[s(t)-\tilde
{s}(t)]|\\
& \leq|u_{x}(x-s^{\prime},t)|t\left\|  v-\tilde{v}\right\|
\end{align*}
where $s^{\prime}$ is the intermediate value between $s(t)$\ \ and $\tilde{s}(t).$

The estimates for the first term in (\ref{u-util}) are based on the maximum
principle. After the free interfaces are determined, both $u$\ \ and
$\tilde{u}$\ \ solve the heat equation off their respective boundaries. Their
difference $w=u-\tilde{u}$ satisfies the heat equation in each of the three
domains
\[
x<s(\tau),\;s(\tau)<x<\tilde{s}(\tau),\;\tilde{s}(\tau)<x,\;\tau<\sigma
\]
here we assumed that $s(\tau)<\tilde{s}(\tau)$ and that $\sigma$ is such that
the inequality holds for all $\tau<\sigma.$ \ It is easy to estimate the
boundary values of the difference
\begin{align*}
\left|  w((s(t),t)\right|   & =\left|  u(s(t),t)-\tilde{u}(s(t),t)\right| \\
& =\left|  g^{-1}(v(t))-g^{-1}(\tilde{v}(t))-\tilde{u}_{x}(s^{\prime
},t)[s(t)-\tilde{s}(t)]\right| \\
& \leq L|v(t)-\tilde{v}(t)|+\left|  \tilde{u}_{x}(s^{\prime},t)\right|
t\left\|  v-\tilde{v}\right\|
\end{align*}
where again $s^{\prime}$ is the intermediate value between $s(t)$\ \ and
$\tilde{s}(t).$ Because of the a priori estimate on the derivative
\[
\left|  \tilde{u}_{x}(s^{\prime},t)\right|  \leq\mathcal{M}+|\tilde{u}%
_{0}|\frac{c}{\sqrt{t}}%
\]
(see Theorem \ref{c1alpha}) we obtain
\[
\left|  w((s(t),t)\right|  \leq\left\|  v-\tilde{v}\right\|  (L+\mathcal{M}%
t+c\sqrt{t}\left\|  \tilde{u}_{0}\right\|  ).
\]
A similar estimate holds for the other interface. Of course the initial data
for $w$ is equal to $u_{0}-\tilde{u}_{0}.$ From the maximum principle for each
of the three domains we obtain that%

\[
|u(.,t)-\tilde{u}(.,t)|\leq C_{1}\left\|  v-\tilde{v}\right\|  +C_{2}\left\|
u_{0}-\tilde{u}_{0}\right\|
\]
\end{proof}

\section{Absorbing set and attractor\label{absorb-sec}}

In this section we use the estimates obtained above to establish existence of
bounded absorbing sets and of the attractor which is compact in the weighted
space of continuous functions. In order to establish compactness of the
attractor we need to make use of the heat losses. It can be verified that most
of the estimates and analytical properties of the solutions can be obtained
without the heat losses. The presence of heat losses only improves the
estimates. On the other hand the problem with the heat losses exhibits uniform
exponential decay in time of the contribution of initial data which is
utilized in the proof of compactness of the attractor.

For our purposes we rephrase Theorem \ref{absorb} in the following form:

\begin{proposition}
Let $0\leq\alpha<\alpha_{time}$ (where $\alpha_{time}$ is defined in
(\ref{atime})), then:\newline (i) The semigroup $T_{2}$ is uniformly
exponentially contracting in $C_{\alpha}$:
\[
r_{X}(t)=\sup_{u^{0}\in X}|T_{2}(t)u^{0}|_{\alpha}\leq C\exp(-\kappa
(\gamma,\alpha)t)N,\quad\kappa>0
\]
where $\kappa(\gamma,\alpha)=\gamma-\alpha^{2}-\alpha V^{0}>0$ for any ball
\[
X=\{u\in C_{\alpha};\quad|u|_{\alpha}\leq N\}
\]
\newline (ii) For any $\varepsilon>0,$ the ball $B_{a}:=\{u\in C_{\alpha
}:\quad|u|_{\alpha}\leq\dfrac{V^{0}}{\sqrt{\gamma}}+\varepsilon\}$ is an
absorbing set for bounded subsets of $C_{\alpha}.$ Here the radius $a$ of the
absorbing ball reflects the contribution of the free interface alone.
\end{proposition}

Next we prove that the boundary contribution to the evolution, i.e. the
operators $T_{1}(t)$ are \textit{uniformly compact}. Namely, the following
proposition holds:

\begin{proposition}
If $\alpha<\alpha_{space}=\min(\dfrac{v_{0}}{4},\dfrac{\gamma}{2V^{0}})$ then
for any $t_{0}>0$ the orbit of the ball $\cup_{t\geq t_{0}}T_{1}(t)X$ is
relatively compact in $C_{\alpha}.$
\end{proposition}

\begin{proof}
For \ the version of Arzela-Ascoli theorem appropriate for $C_{\alpha}$ it is
sufficient to have uniform boundedness for the derivative and uniform decay of
the family of functions as $|x^{\prime}|\rightarrow\infty$ which is faster
than the decay prescribed by the weight$.$ From the estimate
(\ref{T1-solution}) we see that the contributions from the interface decay as
$\ \exp(-(\alpha_{space}+\varepsilon)|x|)$ which can be made faster than any
$\exp(-\alpha|x|)$ for $\alpha<\alpha_{space}.$ On the other hand, the
weighted estimate (\ref{T1-derivative}) (cf. also Remark \ref{4+epsilon})
demonstrate that the spatial derivative is uniformly bounded. Then it is easy
to construct a finite $\varepsilon$-net by choosing a finite interval beyond
which the functions of the family are smaller than $\varepsilon$ and extending
the elements of the $\varepsilon$-net from this interval by zero.
\end{proof}

The properties of the evolution operator $T(t)$ described in the above
propositions allow us to apply the abstract general result (see, for example,
\cite[Chap. 1]{temam}) that in our situation can be stated as follows:

\begin{theorem}
The $\omega$-limit set $\mathcal{A}_{\mathcal{\alpha}}$ of the absorbing set
$B_{a}$ is a global exponential compact attractor for the metric space
$C_{\alpha}$; $\mathcal{A}_{\mathcal{\alpha}}$ is the maximal attractor in
$C_{\alpha}$ and it is connected.
\end{theorem}

In order to demonstrate extra regularity of the elements of the attractor, we
shall need a general fact concerning compact attractors (cf. \cite{chueshov}).

\begin{theorem}
Let $X$ be a Banach space and $B\subset X$ be a $\ $ball. Let $S:B\rightarrow
B$ be a uniformly continuous mapping and $\mathcal{C}\subset B $ a compact
attractor for the iteration semigroup $S^{n}.$ Then on the attractor $S$ is a
mapping onto $S\mathcal{C=C}$.
\end{theorem}

The theorem holds for either continuous or discrete time. For the simplicity
of presentation we consider only the discrete case here.

\begin{proof}
Suppose $S\mathcal{C\neq C}$, then there exists $x_{0}\in\mathcal{C}$ and
$x_{0}\notin S\mathcal{C}$. In this case there exists a whole ball
$B_{r}(x_{0})=\{x\in\mathcal{C}:\left\|  x-x_{0}\right\|  <r\}$ such that
$B_{r}(x_{0})\cap S\mathcal{C}=\emptyset.$ Indeed, since the attractor is
compact, its continuous image is compact and therefore the distance to
$x_{0},$ being a continuous function on \ a compact, attains its nonzero
minimum $r$ on $S\mathcal{C}$.

We know that $S$ is uniformly continuous, therefore for any $\varepsilon$
there exists $\delta(\varepsilon)$ so that $\left\|  Sx-Sy\right\|
<\varepsilon$ if $\left\|  x-y\right\|  <\delta.$

Since $x_{0}\in\mathcal{C}$ for any $\varepsilon$ there exists $n$ so that
$\left\|  S^{n}x-x_{0}\right\|  <\varepsilon$ for some $x\in B,$
simultaneously we can always select $n$ so large that $dist(S^{n-1}%
B,\mathcal{C)}<\delta$. Now take $y=S^{n-1}x$ and $a\in\mathcal{C}$, $\left\|
y-a\right\|  <\delta,$ then $\left\|  Sa-x_{0}\right\|  \leq\left\|
Sa-Sy\right\|  +\left\|  Sy-x_{0}\right\|  \leq\varepsilon+\varepsilon.$ By
selecting $\varepsilon<r/2$ we come to a contradiction.
\end{proof}

As an immediate application of the above theorem we obtain the following
regularity result:

\begin{theorem}
\label{a-alpha} For $0\leq\alpha<\alpha_{\min}^{\prime}$ (where $\alpha_{\min
}^{\prime}$\ is defined in (\ref{amin})) the semigroup on the attractor
$\mathcal{A}_{\alpha}$ in the space $C_{\alpha}$ is onto; $\mathcal{A}%
_{\alpha}$ consists of differentiable functions that satisfy the estimates
\[
|\phi|_{\alpha}\leq\dfrac{V^{0}}{\sqrt{\gamma}},\quad|\phi_{x}|_{\alpha}%
\leq\mathcal{M}(v_{0},V^{0},\alpha,\gamma)\leq\mathcal{M}(v_{0},V^{0}%
,\alpha_{\min}^{\prime},\gamma)
\]
which yields that all $\mathcal{A}_{\alpha}\subseteq C_{\beta}^{1}$ for
$0\leq\alpha\leq\beta=\alpha_{\min}^{\prime}-\varepsilon.$
\end{theorem}

\begin{proof}
The theorem follows from the previous one if we note that the required uniform
continuity follows from the well-posedness Theorem \ref{well-pose}. Since the
mapping is onto, given $\phi\in\mathcal{A}_{\alpha}$ for any $t$ there exist
$\psi\in\mathcal{A}_{\alpha}$, so that $\phi=T(t)\psi.$ By using estimates
(\ref{T1-solution})-(\ref{T1-derivative}) \ and taking into account
exponential decay of the contribution from initial data as $t\rightarrow
\infty$ (\ref{T2-solution})-(\ref{T2-derivative}) we obtain the result.
\end{proof}

\begin{remark}
Since any function in the attractor can be viewed as a result of evolution by
the semigroup, it therefore locally satisfies the heat equation and
consequently it is locally $C^{\infty}$. In addition we can show that due to
the differentiability of the velocity of the interface, functions in the
attractor are $C_{\alpha}^{3}$ up to the interface.
\end{remark}

In addition to being onto, the semigroup $T(t)$ is also one-to-one on the attractor:

\begin{proposition}
The semigroup $T(t)$ is one-to-one on the attractor $\mathcal{A}_{\alpha}.$
\end{proposition}

\begin{proof}
Let $T(t)u_{1}=T(t)u_{2}.$ Denote $t_{0}=\inf\{t:(T(t)u_{1})(x)\equiv
(T(t)u_{2})(x)\};$ obviously $t_{0}>0.$ Then the difference $w=T(t)u_{1}%
-T(t)u_{2}$ is identically $0$ for $t\geq t_{0}.$ Let $t_{1}<t_{0} $ then
there exists $x_{0}$ such that $w(x_{0})\neq0.$ next we select the parabolic
neighborhood (a cup) $U(x_{0})=\{(x,t):t_{1}\leq t<t_{0}+\varepsilon
,\;x_{0}-\delta-k(t-t_{1})<x<x_{0}-\delta+k(t-t_{1})\}$ where $\varepsilon,$
and $\delta$ and $k$ are selected in such a way that the neighborhood does not
intersect the free interface. Since $\ w$ is a solution of the heat equation
in $U(x_{0})$ it is real analytic in $t$ and therefore should be identically
$0$ in $U(x_{0})$ in contradiction with the assumption $w(x_{0})\neq0.$
\end{proof}

Computations for the volume evolution estimates and the Hausdorff dimension
below are implemented in a Hilbert space. It is easy to see that the
exponential decay implies the inclusion $C_{\alpha}\subset H_{0}$ see
(\ref{beta-alpha}). Consequently, the compactness result holds for $H_{0}$ as well:

\begin{theorem}
Let $0<\alpha<\alpha_{\min}$. Then:\newline The semigroup $T_{2}$ is uniformly
exponentially contracting in the $H_{0}$-norm:
\[
r_{X}(t)=\sup_{u^{0}\in X}\left\|  T_{2}(t)u^{0}\right\|  _{H_{0}}\leq
C\exp(-\kappa(\gamma,\alpha)t)N/\sqrt{\alpha},\quad\kappa>0
\]
where $\kappa(\gamma,\alpha)=\gamma-\alpha^{2}-\alpha V^{0}>0$ for any ball
\[
X=\{u\in C_{\alpha};\quad|u|_{\alpha}\leq N\}.
\]
\newline For any $t_{0}>0$ the orbit of the ball $\cup_{t\geq t_{0}}T_{1}(t)X$
is relatively compact in $H_{0}.$\newline The $\omega$-limit set
$\mathcal{A}_{\mathcal{\alpha}}$ of the absorbing set $B_{a}:=\{u\in
C_{\alpha}:\quad|u|_{\alpha}\leq\dfrac{V^{0}}{\sqrt{\gamma}}+\varepsilon\}$ is
a relatively compact set in the $H_{0}$-metric.
\end{theorem}

\begin{proof}
The only additional ingredient of the proof, as compared to the $C_{\alpha}$
case is provided by the imbedding estimate (\ref{beta-alpha}), which yields
that the $\varepsilon$-net generated for a set in the $C_{\alpha}$-norm is
automatically an $(\varepsilon/\sqrt{\alpha})$-net in the $H_{0}$-norm. Since
the imbedding estimate implies continuity of the imbedding, the above proof
essentially repeats the proof of the fact that a continuous image of a compact
set is compact.
\end{proof}

\section{Evolution of the volume elements on the attractor}

\setcounter{equation}{0} In this section we present the main result of the
paper which is a proof that the Hausdorff dimension of the attractor is finite
(for definiteness we consider $\mathcal{A}_{0},$i.e. $\alpha=0$)$.$ The proof
is based on a study of evolution of the infinitesimal volume along the
trajectories in the attractor. We demonstrate that for sufficiently large $m$
that is defined solely by the physical parameters of the problem the
$m$-dimensional volume decays exponentially. This property combined with the
compactness suggests that the Hausdorff dimension of the attractor for the
solutions of the free boundary problem is no larger than $m$. In the arguments
regarding the Hausdorff dimension of the attractor we follow quite closely the
ideas outlined in \cite{temam}.

First we restate the problem in the coordinate frame attached to the free
interface, $\tilde{x}=x-s(t)$ as follows
\begin{align}
u_{t}  & =u_{xx}+v(t)u_{x}-\gamma u:=F(u),\quad-\infty<x<\infty,\quad
x\neq0\nonumber\\
g(u|_{x=0})  & =v(t),\quad\lbrack\partial u/\partial x]|_{x=0}=v(t),\quad
\label{attach}\\
u(x,0)  & =u^{0}(x).\nonumber
\end{align}
(Tildes have been omitted.)

Let $\{U(.,t),V(t)\}$ be an orbit in the attractor. Let us consider the formal
linearization of the problem (\ref{attach}) about $\{U,V\}$:%

\begin{equation}
z_{t}=z_{xx}+z_{x}V-\gamma w-z(0,t)U_{x}/\nu(V(t)):=F^{\prime}(U,V)z\label{l1}%
\end{equation}%
\begin{equation}
z(0,t)+\nu(V(t))[z_{x}(0,t)]=0,\label{l2}%
\end{equation}%
\begin{equation}
z(x,0)=z_{0}(x)\label{l3}%
\end{equation}
where $\nu(V)=-(g^{-1})^{\prime}(V)$. We require $\nu(V)$ to be positive and
bounded from below; $\nu(V)\geq\nu_{0}$. This condition again mimics the
behavior of the Arrhenius kinetics. We have eliminated the velocity
perturbation $v(t)$ in the term $v(t)U_{x}$ of the \ linearization through
replacing it by the perturbation of the temperature $z(0,t)/\nu(V(t))$ that
arises from the linearization of the kinetic boundary condition in
(\ref{attach}). The linearized problem represents the first variation of
problem (\ref{attach}).

It is possible to show that the linearized problem is well-posed in the
following sense:

\begin{theorem}
\label{linear}For any $z_{0}\in H$ there exists a unique solution $z$ of
(\ref{l1}-\ref{l3}) such that $z\in L^{2}(0,T;\Xi(t))\cap C([0,T];H)$ where
$\Xi(t)=\{f\in H^{1},\ f(0)+[f_{x}(0)]\nu(V(t))=0\}$
\end{theorem}

\begin{proof}
This linear problem is somewhat nonstandard as it contains a nonlocal term
(projection) $z(0,t).$ Nonetheless it can be handled as follows. Consider
first the problem (\ref{l1}-\ref{l3}) with a source, and zero initial
conditions
\[
\widetilde{w}_{t}=\widetilde{w}_{xx}+\widetilde{w}_{x}V-\gamma\widetilde
{w}+\mathcal{F}(x,t)
\]%
\[
\widetilde{w}(0,t)+[\widetilde{w}_{x}(0,t)]\nu(V(t))=0,
\]%
\[
\widetilde{w}(x,0)=0,
\]
and let $\frak{L}$ be its solution operator: $\widetilde{w}=\frak{L}%
\mathcal{F}(x,t)$. Existence of unique global solutions for such problems is
guaranteed by the general theory of linear parabolic equations.

We regard a solution of (\ref{l1}-\ref{l3}) as a superposition of an
appropriate $\widetilde{w}$ and of $W(x,t)$ which solves the
\textit{homogeneous} problem with the initial condition $z_{0}(x)$.Then, with
the nonlocal term viewed as a source, on the boundary one obtains an equation
for $z(0,t)$:
\begin{equation}
\frak{L}[-z(0,t)U_{x}/\nu(V(t))+W(x,t)]_{x=0}=z(0,t)+W(0,t).
\end{equation}
It is not difficult to show that the above equation is uniquely solvable as an
integral equation with a sufficiently regular kernel. Thus, the source term is
found and, consequently the problem (\ref{l1}-\ref{l3}) can be solved.
\end{proof}

In order to estimate the evolution of the volume element we need an estimate
for $||U_{x}||_{H_{0}}.$ From now on, for brevity the $H_{0}$-norm will be
denoted by $||.||.$ From the imbedding estimate (\ref{beta-alpha}) and
estimate (\ref{T1-derivative}) we obtain the following important result:

\begin{lemma}
\label{uxlemma} If $U\in\mathcal{A}_{0},$ then $||U_{x}||\leq\mathcal{M}%
/\sqrt{\alpha_{\min}^{\prime}}:=\mathcal{N},$ where $\alpha_{\min}^{\prime
}=\min(\frac{v_{0}}{8},\frac{\gamma}{2V^{0}})$
\end{lemma}

We are now ready to estimate the evolution of the volume element. To this end
we need to estimate the trace of the finite-dimensional projections of the
generator of the linear semigroup. Let $\{\xi_{1},\ldots,\xi_{m}\}$ be $m$
elements of $H$ and let $\{z_{1},\ldots,z_{m}\}$ be the corresponding
solutions of the linearized problem. Then it can be shown that the volume
element spanned by $\{\xi_{1},\ldots,\xi_{m}\}$ evolves accordingly to the
formula
\[
|z_{1}(t)\wedge\ldots\wedge z_{m}(t)|=|\xi_{1}\wedge\ldots\wedge\xi_{m}%
)|\exp\int\nolimits_{0}^{t}Tr\,[F^{\prime}(U(\tau),V(\tau))\circ Q_{m}%
(\tau)]d\tau,
\]
where $Q_{m}(\tau)=Q_{m}(\tau,U,V;\xi_{1},\ldots,\xi_{m})$ is the projector in
$H$ onto the space spanned by $\Xi(\tau)=\{z_{1}(\tau),\ldots,z_{m}(\tau)\}$.
In order to calculate the trace we need to choose a basis in $\Xi(\tau)$
orthogonal in the sense of $H$.

Evaluation of the inner product in $H$ gives rise to sums of integrals over
the domain $(-\infty,0)\cup(0,\infty).$ For brevity everywhere in the sequel
we denote them by
\begin{equation}
\int_{R^{\pm}}f(x)dx=(\int_{-\infty}^{0}+\int_{0}^{\infty}%
)f(x)dx.\label{intinf}%
\end{equation}

Let $\phi$ be an element of $\Xi(\tau)$. Consider the following inner product
in $H$%

\begin{equation}
\langle F^{\prime}\phi,\phi\rangle=-\gamma\langle\phi,\phi\rangle+\int
_{R^{\pm}}\phi_{xx}\phi dx+V\int_{R^{\pm}}\phi_{x}\phi dx+[\phi^{\prime
}(0)]\int_{R^{\pm}}U_{x}\phi\,dx=-\gamma+I_{1}+I_{2}+I_{3}%
\end{equation}
We integrate $I_{1}$ by parts,
\[
I_{1}=\phi_{x}\phi|_{-\infty}^{0}+\phi_{x}\phi|_{0}^{\infty}-\int_{R^{\pm}%
}\phi_{x}^{2}dx=-[\phi^{\prime}(0)]\phi(0)-\int_{R^{\pm}}\phi_{x}^{2}dx
\]
It is easily seen that%

\[
\int_{R^{\pm}}\phi_{x}\phi dx=\frac{1}{2}\int_{R^{\pm}}(\phi^{2})^{\prime}dx=0
\]
then
\begin{equation}
\langle F^{\prime}\phi,\phi\rangle=-\gamma-[\phi^{\prime}(0)]\phi
(0)-\int_{R^{\pm}}\phi_{x}^{2}dx+[\phi^{\prime}(0)]\int_{R^{\pm}}U_{x}\phi\,dx
\end{equation}

In the choice of the $m$-dimensional orthogonal set we shall distinguish the
two possibilities: $\phi(0)=0$ (which defines an ($m-1)$-dimensional
subspace), and otherwise. Since the trace of the operator is independent of
the choice of an orthonormal basis, we can choose $m-1$ basis elements
satisfying the above condition. In the case $\phi(0)=0$ we obtain
\begin{equation}
\langle F^{\prime}\phi,\phi\rangle=-\gamma-\int_{R^{\pm}}\phi_{x}^{2}%
dx\leq-\gamma\label{elemV-3}%
\end{equation}
Note that the terms with $[\phi^{\prime}(0)]$ vanish since $[\phi^{\prime
}(0)]=-\phi(0)/\nu=0$ in view of the boundary condition.

For the basis element with $\phi(0)\neq0$ the corresponding trace component,
\[
\langle F^{\prime}\phi,\phi\rangle=-\gamma-[\phi_{x}]\phi|_{0}-\int_{R^{\pm}%
}\phi_{x}^{2}dx-\frac{\phi(0)}{\nu}\int_{R^{\pm}}U_{x}\phi dx
\]
is estimated from above as follows. First we estimate the last term:
\[
\left|  \frac{\phi(0)}{\nu}\int_{R^{\pm}}U_{x}\phi dx\right|  \leq\frac{1}%
{2a}\int_{R^{\pm}}U_{x}^{2}dx+\frac{a}{2}\frac{\phi^{2}(0)}{\nu^{2}}%
\int_{R^{\pm}}\phi^{2}dx=\frac{1}{2a}\int_{R^{\pm}}U_{x}^{2}dx+\frac{a}%
{2}\frac{\phi^{2}(0)}{\nu^{2}}%
\]
where $a>0$ will be chosen later on.

Now we need the following interpolation result. By integrating from $-\infty$
to 0 we obtain:
\[
\phi^{2}(0)=\int_{-\infty}^{0}(\phi^{2})_{x}dx=2\int_{-\infty}^{0}\phi\phi
_{x}dx\leq2\int_{-\infty}^{0}(\frac{c\phi_{x}^{2}}{2}+\frac{\phi^{2}}{2c})dx
\]
On the other hand,
\[
\phi^{2}(0)=-\int_{0}^{\infty}(\phi^{2})_{x}dx=-2\int_{0}^{\infty}\phi\phi
_{x}dx\leq2\int_{0}^{\infty}(\frac{c\phi_{x}^{2}}{2}+\frac{\phi^{2}}{2c})dx.
\]
Thus
\[
\phi^{2}(0)\leq\int_{R^{\pm}}(\frac{c\phi_{x}^{2}}{2}+\frac{\phi^{2}}%
{2c})dx=\int_{R^{\pm}}(\frac{c\phi_{x}^{2}}{2})dx+\frac{1}{2c}%
\]
In some sense, this estimate is a Sobolev trace theorem.

Finally we obtain
\begin{align*}
\langle F^{\prime}\phi,\phi\rangle & \leq-\gamma+\frac{\phi^{2}(0)}{\nu}%
-\int_{R^{\pm}}\phi_{x}^{2}dx+\frac{1}{2a}\int_{R^{\pm}}U_{x}^{2}dx+\frac
{a}{2}\frac{\phi^{2}(0)}{\nu^{2}}\\
& \leq-\gamma+(\frac{1}{\nu}+\frac{a}{2\nu^{2}})\frac{1}{2}[c\int_{R^{\pm}%
}\phi_{x}^{2}dx+\frac{1}{c}]-\int_{R^{\pm}}\phi_{x}^{2}dx+\frac{1}{2a}%
\int_{R^{\pm}}U_{x}^{2}dx\\
& =-\gamma+[(\frac{1}{\nu_{0}}+\frac{a}{2\nu_{0}^{2}})\frac{1}{2}%
c-1]\int_{R^{\pm}}\phi_{x}^{2}dx+\frac{1}{2a}\int_{R^{\pm}}U_{x}^{2}%
dx+(\frac{1}{\nu_{0}}+\frac{a}{2\nu_{0}^{2}})\frac{1}{2c}%
\end{align*}
If $a$\ and\ $c$ are chosen so that the coefficient at the integral of
$\phi_{x}^{2}$ is nonpositive, say if $c=4\nu_{0}^{2}/(2\nu_{0}+a)$ then
\[
\langle F^{\prime}\phi,\phi\rangle\leq-\gamma+\frac{1}{2a}\int_{R^{\pm}}%
U_{x}^{2}dx+(\frac{2\nu_{0}+a}{4\nu_{0}^{2}})^{2}%
\]
for any $a>0.$

By Lemma \ref{uxlemma} the norm \ $||U_{x}||\leq\mathcal{N}$ where the bound
depends only on the kinetics. To optimize the estimate above we choose $a$
that gives the minimum to the expression $(\frac{2\nu+a}{4\nu^{2}})^{2}%
+\frac{1}{2a}\mathcal{N}^{2}$ considered as a function of $a.$ This results in
the estimate
\begin{equation}
\langle F^{\prime}\phi,\phi\rangle\leq\mu=-\gamma+\min_{a>0}\left(
(\frac{2\nu_{0}+a}{4\nu_{0}^{2}})^{2}+\frac{1}{2a}\mathcal{N}^{2}\right)
\leq(\frac{2\nu_{0}+1}{4\nu_{0}^{2}})^{2}+\frac{1}{2}\mathcal{N}%
^{2}\label{elemV-9}%
\end{equation}
The explicit form of $\mu$ is not important for our purposes. Thus, employing
the above estimates for the trace entries (\ref{elemV-3}), (\ref{elemV-9}) we
can complete the estimate for the evolution of the volume element:
\begin{equation}
Tr\,[F^{\prime}(U(\tau),V(\tau))\circ Q_{m}(\tau)]=\sum_{i=1}^{m}\langle
F^{\prime}\phi_{i},\phi_{i}\rangle\leq\mu-m\gamma\label{dim2}%
\end{equation}
Taking $m>M=\mu/\gamma$ is sufficient for the trace to become negative. Note
that $M$ depends on $\nu_{0},$ $\gamma,$ $V_{0}$ and $v_{0}$.

\section{Differentiability of the semigroup}

\setcounter{equation}{0} To utilize the trace estimate developed in the
previous section we need to demonstrate that the \textit{nonlinear} evolution
of the volume is well approximated by its \textit{linear} counterpart. This
will be ensured by the differentiability of the semigroup solving the
free-interface problem with respect to the initial conditions, see
\cite{temam}.

For the purposes of this section we need to impose an additional condition on
the kinetics function: we will require both $g$ and $g^{-1}$ to be twice
differentiable. In applications this condition is definitely satisfied for all
realistic kinetics.

In Sec.~2 we cited the global existence result for the \textit{classical}
solutions of the free interface problem (\ref{he})-(\ref{jc}). However, our
trace estimates take place in the geometry of a Hilbert space. Therefore we
need to introduce weak solutions by extending the existence theory to more
general initial data that belong to a Hilbert space. The scheme of
introduction of weak solutions is based on the following estimates:

\begin{proposition}
Let $U$ and $W$ be two orbits (i.e., two solutions of the problem
(\ref{attach})) with initial data $U_{0},$ $W_{0}$ in the attractor:
$U=T(t)U_{0}$, $W=T(t)W_{0}$. Then for any $t>0$
\begin{equation}
\left\|  U(t)-W(t)\right\|  \leq e^{Ct}\left\|  U_{0}-W_{0}\right\|
\label{u-w}%
\end{equation}%
\[
\int\limits_{0}^{t}\left\|  U(t)-W(t)\right\|  _{1}^{2}d\tau\leq
e^{Ct}\left\|  U_{0}-W_{0}\right\|  ^{2}%
\]
where $C$ is a uniform constant$.$
\end{proposition}

\begin{remark}
The above proposition allows us to obtain weak solutions with initial data in
the closure of $\mathcal{A}$ in the $H$-norm. Namely, in a standard fashion we
select a Cauchy sequence of initial conditions in $\mathcal{A}$ and define the
solution as the corresponding limit of smooth solutions. We note that it is
not necessary to take initial data from $\mathcal{A}$; similarly it is
possible to define a weak solution for the initial data in the closure of a
ball in $C_{\alpha}.$\ 
\end{remark}

\begin{proof}
Let $U(x,t)$ and $W(x,t)$ be two solutions of the free boundary problem (in
the frame attached to the free boundary)
\[
U_{t}=U_{xx}+[U_{x}(0,t)]U_{x}-\gamma U,\quad g(U(0,t))=[U_{x}(0,t)],\quad
U(x,0)=U_{0}(x),
\]%
\[
W_{t}=W_{xx}+[W_{x}(0,t)]W_{x}-\gamma W,\quad g(W(0,t))=[W_{x}(0,t)],\quad
W(x,0)=W_{0}(x).
\]
The difference $w=U-W$ solves the following problem
\[
w_{t}=w_{xx}+[U_{x}(0,t)]w_{x}+[w_{x}(0,t)]W_{x}-\gamma w,\quad
\]%
\[
-[w_{x}(0,t)]=g(W(0,t))-g(U(0,t))=-(g(\theta))^{\prime}w(0,t).
\]%
\[
w(x,0)=U_{0}(x)-W_{0}(x).
\]
We also observe that $-(g(\theta))^{\prime}\leq const=C$ while $U$, $W$ and
their $x$-derivatives are uniformly bounded on the attractor.

We multiply the equation throughout by $w$ and integrate to obtain the
following energy estimate for the $H$ norm:
\begin{gather}
\frac{1}{2}\frac{d}{dt}\left\|  w\right\|  ^{2}=\int_{R^{\pm}}w_{xx}%
wdx+[U_{x}(0,t)]\int_{R^{\pm}}w_{x}wdx+[w_{x}(0,t)]\int_{R^{\pm}}%
W_{x}wdx-\gamma\int_{R^{\pm}}w^{2}dx\nonumber\\
=-[w_{x}]w|_{0}-\gamma\left\|  w\right\|  ^{2}-\left\|  w_{x}\right\|
^{2}+[U_{x}(0,t)]\int_{R^{\pm}}w_{x}wdx+[w_{x}(0,t)]\int_{R^{\pm}}%
W_{x}wdx\label{diff3}%
\end{gather}
We need to estimate different terms in (\ref{diff3}). For the first term we
get on respective intervals
\[
|[w_{x}(0,t)]w(0,t)|\leq Cw(0,t)^{2}=C\int\limits_{-\infty}^{0}(w^{2}%
)_{x}dx\leq2C\int\limits_{-\infty}^{0}|w_{x}w|dx
\]%
\[
|[w_{x}(0,t)]w(0,t)|\leq Cw(0,t)^{2}=-C\int\limits_{0}^{\infty}(w^{2}%
)_{x}dx\leq2C\int\limits_{0}^{\infty}|w_{x}w|dx
\]
The sum of the above inequalities yields the estimate
\begin{equation}
|[w_{x}(0,t)]w(0,t)|\leq Cw(0,t)^{2}\leq C(\varepsilon_{1}\left\|  w\right\|
^{2}+\frac{1}{\varepsilon_{1}}\left\|  w_{x}\right\|  ^{2})\label{poincare}%
\end{equation}

Next,
\[
|[U_{x}(0,t)])\int_{R^{\pm}}w_{x}wdx|\leq|[U_{x}(0,t)]|(\varepsilon
_{2}\left\|  w\right\|  ^{2}+\frac{1}{\varepsilon_{2}}\left\|  w_{x}\right\|
^{2})\leq C_{1}(\varepsilon_{2}\left\|  w\right\|  ^{2}+\frac{1}%
{\varepsilon_{2}}\left\|  w_{x}\right\|  ^{2})
\]
Also,
\begin{gather*}
\left|  \lbrack w_{x}(0,t)]\int_{R^{\pm}}W_{x}wdx\right|  \leq C_{3}%
|w(0,t)|\ \left\|  W_{x}\right\|  \ \left\|  w\right\| \\
\leq C_{3}|w(0,t)|\ \mathcal{N}\ \left\|  w\right\|  \leq C_{4}(\varepsilon
_{3}\left\|  w\right\|  ^{2}+\frac{1}{\varepsilon_{3}}\left\|  w_{x}\right\|
^{2})^{1/2}\left\|  w\right\| \\
\leq C_{4}(\sqrt{\varepsilon_{3}}\left\|  w\right\|  +\frac{1}{\sqrt
{\varepsilon_{3}}}\left\|  w_{x}\right\|  )\left\|  w\right\|  =C_{4}%
(\sqrt{\varepsilon_{3}}\left\|  w\right\|  ^{2}+\frac{1}{\sqrt{\varepsilon
_{3}}}\left\|  w_{x}\right\|  \left\|  w\right\|  )\\
\leq C_{4}(\varepsilon_{5}\left\|  w\right\|  ^{2}+\frac{1}{\varepsilon_{4}%
}\left\|  w_{x}\right\|  ^{2})
\end{gather*}

Collecting the estimates for different terms we get
\begin{gather*}
\frac{1}{2}\frac{d}{dt}|w|_{\alpha}^{2}\leq-\left\|  w_{x}\right\|
^{2}-\gamma\left\|  w\right\|  ^{2}\\
+C(\varepsilon_{1}\left\|  w\right\|  ^{2}+\frac{1}{\varepsilon_{1}}\left\|
w_{x}\right\|  ^{2})+C_{1}(\varepsilon_{2}\left\|  w\right\|  ^{2}+\frac
{1}{\varepsilon_{2}}\left\|  w_{x}\right\|  ^{2})+C_{5}(\varepsilon
_{5}\left\|  w\right\|  ^{2}+\frac{1}{\varepsilon_{4}}\left\|  w_{x}\right\|
^{2})\\
\leq-\frac{1}{2}\left\|  w_{x}\right\|  ^{2}+C_{0}\left\|  w\right\|  ^{2}%
\end{gather*}
where on the last step we have chosen $C/\varepsilon_{1}+C_{1}/\varepsilon
_{2}+C_{5}/\varepsilon_{4}<1/2$. We rewrite our last result as
\begin{equation}
\frac{d}{dt}\left\|  w\right\|  ^{2}+\left\|  w_{x}\right\|  ^{2}\leq
C\left\|  w\right\|  ^{2}.\label{gron1}%
\end{equation}
From this inequality we get first that $\dfrac{d}{dt}\left\|  w\right\|
^{2}\leq C\left\|  w\right\|  ^{2}$ yielding, by Gronwall's inequality, that
\[
\left\|  w\right\|  ^{2}\leq\left\|  w_{0}\right\|  ^{2}\exp(Ct);
\]
at the same time by rearranging and integrating (\ref{gron1}) we obtain
\[
\int\limits_{0}^{t}\left\|  w_{x}\right\|  ^{2}d\tau\leq\int\limits_{0}%
^{t}(C\left\|  w\right\|  ^{2}-\frac{d}{dt}\left\|  w\right\|  ^{2})d\tau
\leq\left\|  w_{0}\right\|  ^{2}\exp(Ct).
\]
\end{proof}

If the initial data are in $H^{1}\cap C_{\alpha}$ then a similar argument
yields the following estimate analogous to (\ref{u-w}):
\begin{equation}
\left\|  U(t)-W(t)\right\|  _{1}\leq e^{Ct}\ \left\|  U_{0}-W_{0}\right\|
_{1}\label{wh1}%
\end{equation}

Next we prove the differentiability in $H^{1}$ that is sufficient for the
validity of the dimension estimate because it implies the differentiability on
$\mathcal{A}\subset H^{1}.$

\begin{theorem}
Let $U$ and $W$ be two orbits $U=T(t)U_{0}$, $W=T(t)W_{0}$,\ \ $U_{0},W_{0}\in
H^{1}\cap C_{\alpha}$. Then there exists $z(t)$ such that
\[
\left\|  U(t)-W(t)-z(t)\right\|  \leq const\ \left\|  U_{0}-W_{0}\right\|
_{1}^{2}%
\]
as $W_{0}\rightarrow U_{0}$.
\end{theorem}

In this case the Frech\'{e}t differential of $T(t)$ at the point $U_{0}$ is
the mapping $z(0)=U_{0}-W_{0}\rightarrow z(t)$, where $z(t)$ solves the
linearized problem.

\begin{proof}
The goal of the proof is to evaluate the difference between $w=U-W$ and its
approximation by the differential. We define $z(x,t)$ as a solution of the
free-interface problem linearized about the orbit $U(x,t)$:
\begin{align}
z_{t}  & =z_{xx}+[z_{x}(0,t)]U_{x}+[U_{x}(0,t)]z_{x},\quad\label{lin}\\
z(0,t)  & =(g^{-1})^{\prime}([U_{x}(0,t)])[z_{x}(0,t)],\quad z(x,0)=U_{0}%
(x)-W_{0}(x),\nonumber
\end{align}
(see Theorem \ref{linear}). For the difference $y=w-z$ we have the following
equations
\begin{align*}
y_{t}  & =y_{xx}+[y_{x}(0,t)]U_{x}+[U_{x}(0,t)]y_{x},\quad\\
y(0,t)  & =(g^{-1})^{\prime}([U_{x}(0,t)])[y_{x}(0,t)]+(g^{-1})^{\prime\prime
}(\theta)[w_{x}(0,t)]^{2}/2,\quad y(x,0)=0,
\end{align*}
We multiply the equation throughout by $y$ and integrate to obtain the
following identity for the $H$ norm:
\begin{align}
\frac{1}{2}\frac{d}{dt}\left\|  y\right\|  ^{2}  & =\int_{R^{\pm}}%
y_{xx}ydx+[U_{x}(0,t)]\int_{R^{\pm}}y_{x}ydx+[y_{x}(0,t)]\int_{R^{\pm}}%
U_{x}ydx\label{y1}\\
& =-[y_{x}]y|_{0}-\left\|  y_{x}\right\|  ^{2}+[U_{x}(0,t)]\int_{R^{\pm}}%
y_{x}ydx+[y_{x}(0,t)]\int_{R^{\pm}}U_{x}ydx\nonumber
\end{align}
We need to estimate different terms in (\ref{y1})
\begin{gather*}
|[y_{x}(0,t)]y(0,t)|\leq Cy^{2}(0,t)+C[w_{x}(0,t)]^{2}|y(0,t)|\leq\frac{3}%
{2}Cy^{2}(0,t)+\frac{1}{2}C[w_{x}(0,t)]^{4}\\
=B_{1}y^{2}(0,t)+B_{2}[w_{x}(0,t)]^{4}\\
\leq\dfrac{B_{1}}{2}\left|  \int\limits_{-\infty}^{0}y_{x}^{2}dx\right|
+\dfrac{B_{1}}{2}\left|  \int\limits_{0}^{\infty}y_{x}^{2}dx\right|
+B_{3}(\varepsilon_{1}\left\|  w\right\|  ^{2}+\frac{1}{\varepsilon_{1}%
}\left\|  w_{x}\right\|  ^{2})^{2}\\
\leq B_{1}(\varepsilon_{1}\left\|  y\right\|  ^{2}+\frac{1}{\varepsilon_{1}%
}\left\|  y_{x}\right\|  ^{2})+B_{3}(\varepsilon_{1}\left\|  w\right\|
^{2}+\frac{1}{\varepsilon_{1}}\left\|  w_{x}\right\|  ^{2})^{2}%
\end{gather*}
Next,
\[
\left|  \lbrack U_{x}(0,t)]\int_{R^{\pm}}y_{x}ydx\right|  \leq\left|  \lbrack
U_{x}(0,t)]\right|  (\varepsilon_{2}\left\|  y\right\|  ^{2}+\frac
{1}{\varepsilon_{2}}\left\|  y_{x}\right\|  ^{2})\leq C_{1}(\varepsilon
_{2}\left\|  y\right\|  ^{2}+\frac{1}{\varepsilon_{2}}\left\|  y_{x}\right\|
^{2})
\]
Also,
\begin{gather*}
\left|  \lbrack y_{x}(0,t)]\int_{R^{\pm}}U_{x}ydx\right| \\
\leq(C_{3}|y(0,t)|+B_{4}[w_{x}(0,t)]^{2})\int_{R^{\pm}}|U_{x}y|dx\leq
(C_{3}|y(0,t)|+B_{4}[w_{x}(0,t)]^{2})\ \left\|  y\right\|  \left\|
U_{x}\right\| \\
\leq C_{4}(\varepsilon_{3}\left\|  y\right\|  ^{2}+\frac{1}{\varepsilon_{3}%
}\left\|  y_{x}\right\|  ^{2})^{1/2}\left\|  y\right\|  +B_{5}(\varepsilon
_{1}\left\|  w\right\|  ^{2}+\frac{1}{\varepsilon_{1}}\left\|  w_{x}\right\|
^{2}+\left\|  w\right\|  ^{2})\left\|  y\right\| \\
\leq C_{4}(\sqrt{\varepsilon_{3}}\left\|  y\right\|  +\frac{1}{\sqrt
{\varepsilon_{3}}}\left\|  y_{x}\right\|  )\left\|  y\right\|  +B_{5}%
(\varepsilon_{4}\left\|  w\right\|  ^{2}+\frac{1}{\varepsilon_{1}}\left\|
w_{x}\right\|  ^{2})\left\|  y\right\| \\
\leq C_{5}(\varepsilon_{5}\left\|  y\right\|  ^{2}+\frac{1}{\sqrt
{\varepsilon_{3}}}\left\|  y_{x}\right\|  ^{2})+B_{6}(\left\|  w_{x}\right\|
^{4})
\end{gather*}
where the constants $C_{4}$ and $B_{5}$ include the factor $\left\|
U_{x}\right\|  .$ Collecting the estimates for different terms we get
\begin{gather*}
\frac{1}{2}\frac{d}{dt}\left\|  y\right\|  ^{2}\leq\\
-\left\|  y_{x}\right\|  ^{2}+B_{1}(\varepsilon_{1}\left\|  y\right\|
^{2}+\frac{1}{\varepsilon_{1}}\left\|  y_{x}\right\|  ^{2})+B_{3}%
(\varepsilon_{1}\left\|  w\right\|  ^{2}+\frac{1}{\varepsilon_{1}}\left\|
w_{x}\right\|  ^{2})^{2}\\
+C_{1}(\varepsilon_{2}\left\|  y\right\|  ^{2}+\frac{1}{\varepsilon_{2}%
}\left\|  y_{x}\right\|  ^{2})+C_{5}(\varepsilon_{5}\left\|  y\right\|
^{2}+\frac{1}{\sqrt{\varepsilon_{3}}}\left\|  y_{x}\right\|  ^{2}%
)+B_{6}(\left\|  w_{x}\right\|  ^{4})
\end{gather*}
Now we can select $\varepsilon_{1},$ $\varepsilon_{2},$ and $\varepsilon_{3}$
sufficiently large so that the coefficient by $\left\|  y_{x}\right\|  ^{2}$
is less than $-1/2$. We collect the like terms in the above inequality to
obtain
\[
\frac{1}{2}\frac{d}{dt}\left\|  y\right\|  ^{2}\leq-\frac{1}{2}\left\|
y_{x}\right\|  ^{2}+C\left\|  y\right\|  ^{2}+B\left\|  w\right\|  _{1}^{4}%
\]
We rewrite our last result as
\begin{equation}
\frac{d}{dt}\left\|  y\right\|  ^{2}+\left\|  y_{x}\right\|  ^{2}\leq
C\left\|  y\right\|  ^{2}+B\left\|  w\right\|  _{1}^{4}%
\end{equation}
from where it is clear that
\begin{equation}
\frac{d}{dt}\left\|  y\right\|  ^{2}\leq C\left\|  y\right\|  ^{2}+B\left\|
w\right\|  _{1}^{4}.
\end{equation}
By Gronwall's inequality it yields
\begin{align*}
\left\|  y\right\|  ^{2}  & \leq B\exp(Ct)\int\limits_{0}^{t}\left\|
w\right\|  _{1}^{4}\exp(-C\tau)d\tau\\
& \leq B_{7}\exp(Ct)\left\|  w_{0}\right\|  _{1}^{4}%
\end{align*}
In the above estimate we utilized (\ref{wh1}).
\end{proof}

Finally (see \cite{temam}), the estimate for the dimension of the linear
volume element and differentiability of the semigroup yield the estimate for
the Hausdorff dimension of the attractor:

\begin{theorem}
The Hausdorff dimension of the attractor $\mathcal{A}$ is no larger than
\[
M=[(\frac{2\nu_{0}+1}{4\nu_{0}^{2}})^{2}+\frac{1}{2}\mathcal{N}^{2}]/\gamma
\]
cf. (\ref{dim2}).
\end{theorem}

In conclusion it is worth mentioning that the estimate exhibits a transparent
and physically natural dependence of the dimension on the heat loss and
characteristics of the kinetics which are the defining factors of the
dynamics. We note however that numerical simulations \cite{comp2} on
(\ref{he})-(\ref{jc}) show that the behavior without heat losses and with
sufficiently low heat losses are qualitatively identical and exhibit the same
variety of complex dynamical patterns.

\section{Acknowledgments}

The authors would like to acknowledge support in part by NSF through grants
DMS-0207308 and DMS-9704325.


\begin{thebibliography}{99}
\bibitem{brauner}C.-M. Brauner, J. Hulshof and A. Lunardi, A general approach
to stability in free boundary problems, \textit{J. Differential Equations} 164
(2000), 16--48.

\bibitem{brauner1}C.-M. Brauner and A. Lunardi, Instabilities in a
two-dimensional combustion model with free boundary, \textit{Arch. Ration.
Mech. Anal.} \textbf{154 }(2000)\textbf{,} 157--182.

\bibitem{chen}X. Chen and F. Reitich, Local existence and uniqueness of
solutions of the Stefan problem with surface tension and kinetic undercooling,
\textit{J. Math. Analysis Appl.} 164 (1992) 350-362.

\bibitem{const}P. Constantin, C. Foias, and R. Temam, Attractors representing
turbulent flows, Memoirs of AMS, \textbf{53,} no. 314 (1985)

\bibitem{chueshov}I. D. Chueshov, \textit{Introduction to the Theory of
Infinite-Dimensional Dissipative Systems}, ''Acta'' Scientific Publishing
House (2002).

\bibitem{yin1}S. Din, Y. Tao and H.-M. Yin, A chemical diffusion process with
reaction taking place at free boundary, \textit{Canad. Appl. Math. Quarterly}
5 (1997) 49-74.

\bibitem{minqu}M. Frankel, M. Qu and V. Roytburd, On a free interface problem
modeling solid combustion and rapid solidification in infinite medium,
\textit{World Scientific Series in Applicable Analysis, Vol. 4: Dynamical
Systems and Applications}, R. P. Agarwal ed., World Scientific (1995), 263-278.

\bibitem{siap}M. Frankel, V. Roytburd, and G. Sivashinsky, A sequence of
period doubling and chaotic pulsations in a free boundary problem modeling
thermal instabilities, \textit{SIAM J. Appl. Math.} \textbf{54} (1994), 1101-1112.

\bibitem{port2}M. Frankel V. Roytburd G. Sivashinsky, Complex dynamics
generated by a sharp interface model of self-propagating high-temperature
synthesis, \textit{Combust. Theory Modelling} , 2 (1998), 479-496.

\bibitem{physd}M. Frankel, G. Kova\v{c}i\v{c}, V. Roytburd, and I. Timofeyev,
Finite-dimensional dynamical system modeling thermal instabilities,
\textit{Physica} D 137 (2000), 295-315.

\bibitem{amlet}M. Frankel and V. Roytburd, Finite-dimensional attractors for a
free-boundary problem with a kinetic condition, \textit{\ Appl. Math. Lett.
\textbf{15} (2002), 83-87.}

\bibitem{comp}M. Frankel and V. Roytburd, Compact attractors for a Stefan
problem with kinetics, \textit{EJDE} (2002), 1-27.

\bibitem{sell}M. Frankel and V. Roytburd, Finite-dimensional attractor for a
one-phase Stefan problem with kinetics, \textit{J. Dynam. Differential
Equations \textbf{15 }(2003), 87-106.}

\bibitem{comp2}M. Frankel and V. Roytburd, On attractors for a sharp interface
model of exothermic phase transitions, \textit{Advances Math Sciences
Applications \textbf{14} (2004)\textbf{, }25-40.}

\bibitem{friedman}A. Friedman, \textit{Partial Differential Equations of
Parabolic Type}, Prentice-Hall, Englewood Cliffs, N.J. (1964).

\bibitem{gt}S. M. Gol'berg and M. I. Tribelskii, On laser induced evaporation
of nonlinear absorbing media, \textit{Zh. Tekh. Fiz. (Sov. Phys.-J. Tech.
Phys.)} \textbf{55}, 848-857 (1985).

\bibitem{langer}J. S. Langer, Lectures in the theory of pattern formation,
in:\textit{\ Chance and Matter}, J. Souletie, J. Vannimenus and R. Stora,
eds., Elsevier Science Publishers (1987).

\bibitem{lorenzi}L. Lorenzi, A free boundary problem stemmed from combustion
theory. Part I: Existence, uniqueness and regularity results, \textit{J. Math.
Anal. Appl}. \textbf{274} (2002) 505-535.

\bibitem{luc}S. Luckhaus, Solutions for the two-phase Stefan problem with the
Gibbs--Thomson law for the melting temperature, \textit{Euro. J. of Appl.
Math.} 1 (1990) 101-111.

\bibitem{matsiv-solid}B.J. Matkowsky and G. I. Sivashinsky, Propagation of a
pulsating reaction front in solid fuel combustion, \textit{SIAM J. Appl. Math.
}35 (1978), 230-255.

\bibitem{munir}Z. A. Munir and U. Anselmi-Tamburini, Self-propagating
exothermic reactions: the synthesis of high -temperature materials by
combustion, \textit{Mater. Sci. Rep.} 3 (1989), 277-365.

\bibitem{radkevich}E. V. Radkevich, The Gibbs--Thomson correction and
conditions for the existence of a classical solution of the modified Stefan
problem, \textit{Soviet Math. Dokl.} 41 (1991) 274-278.

\bibitem{shk}K. G. Shkadinsky, B. I. Khaikin and A. G. Merzhanov, Propagation
of a Pulsating Exothermic Reaction Front in the Condensed Phase,
\textit{Combust. Expl. Shock Waves} \textbf{7} (1971), 15-22.

\bibitem{temam}R. Temam, \textit{Infinite-Dimensional Dynamical Systems in
Mechanics and Physics}, Springer Verlag: New York (1988).

\bibitem{saarloos}Van Saarloos, W. and Weeks, J. , \textit{Surface undulations
in explosive crystallization: a nonlinear analysis of a thermal instability},
Physica D \textbf{12 }(1984), 279-294.

\bibitem{var2}A. Varma, A. S. Rogachev, A. S. Mukasyan and S. Hwang,
Combustion synthesis of advanced materials, \textit{Adv. Chem. Eng.} 24
(1998), 79-226.

\bibitem{var1}A. Varma, Form from fire, \textit{Sci. American} Aug. (2000), 58-61.

\bibitem{zel}Ia. B. Zeldovich, G. I. Barenblatt, V. B. Librovich and G. M.
Makhviladze, \textit{The Mathematical Theory of Combustion and Explosions,}
Consultants Bureau, Washington, D.C., 1985.
\end{thebibliography}
\end{document}